\documentclass[12pt,a4paper]{amsart}

\usepackage{amsmath,amsthm,amssymb,stmaryrd,commath}
\usepackage{svg}
\usepackage{pgfplots}
\pgfplotsset{compat=1.17}
\usepackage[font=small,labelfont=bf]{caption}
\usepackage[labelformat=simple]{subcaption}
\usepackage{enumitem}
\usepackage{comment}
\usepackage{cancel}
\usepackage[normalem]{ulem}
\usepackage{graphicx}
\usepackage{xcolor}
\usepackage{csquotes}
\MakeOuterQuote{"}
\usepackage[margin=2.7cm]{geometry}
\usepackage{setspace}
\usepackage{tikz}
\usetikzlibrary{cd, babel}
\usepackage{algorithm}
\usepackage{algpseudocode}

\usepackage[colorlinks]{hyperref}
\usepackage[capitalize]{cleveref}

\hypersetup{
	colorlinks=true,
	linkcolor=blue,
	filecolor=blue,
	citecolor = blue,    
	urlcolor=blue,
}

\newtheorem{thm}{Theorem}[section]
\newtheorem{cor}[thm]{Corollary}
\newtheorem{prop}[thm]{Proposition}
\newtheorem{lem}[thm]{Lemma}

\theoremstyle{definition}
\newtheorem{defn}[thm]{Definition}

\newtheorem{exmp}[thm]{Example}

\theoremstyle{remark}
\newtheorem{rem}[thm]{Remark}

\newcommand{\R}{\mathbb{R}}
\newcommand{\Z}{\mathbb{Z}}
\newcommand{\N}{\mathbb{N}}

\newcommand{\im}{\operatorname{im}}
\newcommand{\F}{\mathbb{F}}
\newcommand{\C}{\operatorname{C}}

\newcommand{\Int}{\operatorname{Int}}

\newcommand{\VR}{\operatorname{VR}}
\newcommand{\Cech}{\operatorname{\check{C}}}
\newcommand{\sd}{\operatorname{sd}}

\newcommand{\birth}{\operatorname{birth}}
\newcommand{\Cyc}{\operatorname{Cyc}}
\newcommand{\T}{\operatorname{T}}
\newcommand{\UTB}{\operatorname{UT}}

\definecolor{nice_blue}{HTML}{70BDFB}

\title[Cycling Signatures]{Cycling Signatures: Identifying Cycling Motions in Time Series using Algebraic Topology}
\author[Ulrich Bauer, David Hien, Oliver Junge and Konstantin Mischaikow]{Ulrich Bauer\textsuperscript{1,2}, David Hien\textsuperscript{1}, Oliver~Junge\textsuperscript{1}\\ and Konstantin~Mischaikow\textsuperscript{3}}

\address{\textsuperscript{1} School of Computation, Information and Technology, Technical University of Munich, 85748 Garching, Germany}
\address{\textsuperscript{2} Munich Data Science Institute, Technical University of Munich, 85748 Garching, Germany}
\address{\textsuperscript{3} Department of Mathematics, Rutgers, The State University of New Jersey, 110 Frelinghusen Rd., Piscataway, NJ 08854-8019, USA}

\begin{document}

\begin{abstract}
	Recurrence is a fundamental characteristic of dynamical systems with complicated behavior.  
	Understanding the inner structure of recurrence is challenging, especially if the system has many degrees of freedom and is subject to noise.
	We develop algebraic topological notions for identifying and classifying elementary recurrent motions -- called \emph{cycling} -- and the transitions between those.
	Statistics on these cycling motions can be computed from sampled trajectories (time series data), providing coarse global information on the structure of the recurrent behavior.
	We demonstrate this through three examples; in particular, we identify and analyze six cycling motions in a four-dimensional system with a hyperchaotic attractor.
	We see this as a promising approach to reveal coarse-grained dynamical information on high-dimensional systems.  
\end{abstract}

\maketitle

\section{Introduction}
Time series data is ubiquitous in applied dynamical systems \cite{KaSch04}.
In this paper, we focus on multivariate time series that exhibit recurrent behavior.
Recurrence is a fundamental characteristic of complex dynamical systems:
typical trajectories repeatedly attain states close to previously visited states.
A natural approach to characterizing this recurrence is to identify trajectories which return
to their initial state, i.e., periodic orbits. 
In fact, simple chaotic dynamics are often described in terms of a scaffolding based on periodic orbits \cite{GuHo-83}.
The advantage of periodic orbits is that, at least abstractly, they are simple objects -- closed trajectories that have the topology of a circle.

The contribution of this paper (see Section~\ref{sec:math} for details) is to employ notions from algebraic topology to automatically and robustly detect and classify time series segments that are almost periodic in a topological sense: A time series segment is called \emph{cycling} if it is `non-trivially' almost-closed and its \emph{cycling signature} captures the types of `cycling' dynamics that the segment contains.

More precisely, to a given time series $\Gamma$, we assign a vector space 
of possible `cycling motions', with the property that for any segment $\gamma$ of $\Gamma$, we get a subspace $\Cyc(\gamma)$ of this vector space which we call the cycling signature. 
Intuitively, the cycling signatures classify the segments of $\Gamma$ according to the holes that $\gamma$ `wraps' around.
The construction of these vector spaces
is based on the degree 1 homology of a suitable cell complex constructed from the data.

The advantage of this approach is that algebraic topology provides well established, robust, dimension independent tools to identify and quantify the required `holes' in the time series data \cite{Ghrist.2007,edelsbrunner2008,carlsson2009,MR2121296}.
In particular, they allow us to compute degree 1 homology (with coefficients in a field), i.e., a vector space whose dimension can loosely be interpreted as `the number of holes' in the data set.
While the dimension of a vector space equals the number of basis elements, in general, there is no unique or canonical choice of basis.
This is the case in our context:
given a segment $\gamma$ of the time series $\Gamma$, the dimension of $\Cyc(\gamma)$, called the \emph{cycling rank} of~$\gamma$, indicates the number of linearly independent 
`holes' that $\gamma$ wraps around, but an identification of the specific `holes' would still depend on a choice of basis.

Insight into the structure of recurrence of a dynamical system can be gained, without making any choice of basis, by analyzing $\Cyc(\gamma)$ of a large number of segments $\gamma$.
For example, the distribution of cycling ranks allows us to infer the number of holes that typical segments `wrap' around.
Segments with cycling rank 1 are of particular interest since two segments with the same 1 dimensional signature wrap around the same holes.
In particular, any typical rank 1 signature corresponds to a type of cycling motion in the data set $\Gamma$.
Higher rank signatures can be used to analyze the relations of and transition between signatures of lower rank. 

Even though the analysis of nonlinear time series data has a long history~\cite{KaSch04,Bradley2015,Marwan2007,MR2787567,de_Silva_2009}, 
to the best of our knowledge, cycling signatures are a complementary tool with two distinguishing features:
First, they can be applied to general time series data without the need for visualizations or prior analysis, e.g., a Poincare section, identification of dominant frequencies, etc. 
Second, the characterization is global in the sense that it detects the topologically different cycling motions and the transitions between them. In this sense, the information that we obtain provides information 
on the internal structure of some recurrent invariant set of the underlying (unknown) dynamical system.
The resulting statistics can be interpreted without detailed knowledge of algebraic topology.
In particular, the computational results in Section~\ref{sec:application} can be interpreted with the informal description provided thus far.

\section{Cycling segments and their classification}
\label{sec:math}

Consider a time-continuous dynamical system given by a smooth flow $\Phi\colon \R\times X\to X$ on a state space $X\subset\R^d$.  Assume that we do not know $\Phi$ but only have access to some (possibly perturbed) sample from this system.  This means that we are given a time series $\Gamma = (x_i)_{i=1,\dots,N}$ together with sample times $T = (t_i)_{i=1,\dots,N}$ such that
\[
	x_i \approx \Phi(t_i-t_{i-1}, x_{i-1}), \qquad i=1,\ldots,N,
\]
where $x_0\in X$ is an initial point, the sample times are assumed to satisfy $t_0 < t_1 < \cdots < t_N$, and $x_i$ are measurements from the underlying flow that may be perturbed by e.g. measurement errors or numerical integration.
A \emph{segment} $\gamma$ of $\Gamma$ is a consecutive subsequence $ x_k,\, x_{k+1},\dots,\,x_\ell$ for some $ k< \ell $.
The \emph{time span} of a segment is $\tau(\gamma)=t_\ell - t_k$.

If $\Gamma$ is sampled from a recurrent set with complicated dynamics, then, given a sufficiently long time span and sampling density of $\Gamma$, the time series will closely return to certain states, i.e., there are tuples $ (i,j) $ such that $ x_{i} \approx x_{j} $, but also points $x_\ell$ with $i<\ell<j$ far away from both $x_{i}$ and $x_{j}$. 
We will later define certain segments containing such a 'recurrence point' $x_i$ as \emph{cycling}. 
A time series may have an abundance of cycling segments. 
However, often these cycling segments can be structured into a small number of different qualitative classes.
We first illustrate these concepts using the classical Lorenz system before we show how homology can be used to identify cycling segments and classify them.

\begin{exmp}[Lorenz system]\label{ex:Lorenz}
    \begin{figure}[h]
    	\centering
    	\includegraphics[width=.8\linewidth]{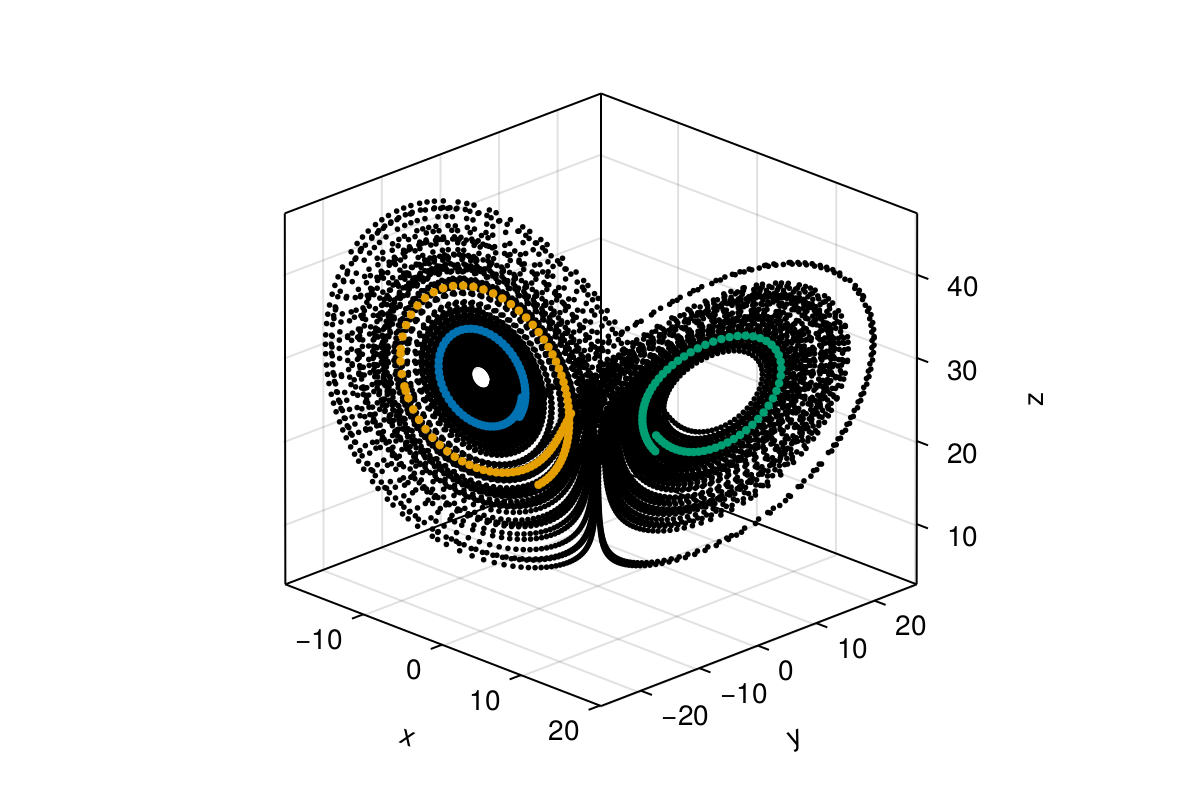}
    	\caption{Cycling segments in the classical Lorenz attractor. 
        The black points come from numerical integration of the system over a large time span. 
        Three segments are highlighted.
        The orange and blue segments belong to the same cycling class which should be different from the cycling class of the green segment. 
        }
    	\label{fig:lorenz-motivating-example}
    \end{figure}
    We consider a time series $\Gamma$ sampled from the Lorenz system with classical parameters (see \cref{sec:Lorenz} for details). \Cref{fig:lorenz-motivating-example} shows segments, colored in orange, blue and green, together with $\Gamma$ which is shown in black. 
    
    Each highlighted segment traces an almost-closed path, these are precisely the type of segments that the notion of cycling is supposed to capture. Both the orange and blue segments wrap around the left `hole' of the attractor while the green segment wraps around the right `hole'. 
    Using our concept of cycling spaces, we will detect that the orange segment is qualitatively similar to the blue segment and different from the green segment.
    \qed
\end{exmp}

\subsection{Cycling Segments}

Our definition of cycling is motivated by the observation that homology can be used to identify periodic orbits. Consider a fragment $\zeta = \Phi([a,b],x)\subset X$ of the orbit through some $x\in X$. 
Then $\zeta$ contains a full periodic orbit if and only if it is a circle, which in turn is the case if and only if $H_1(\zeta)\neq 0$. This holds for homology with coefficients in an arbitrary ring, and from now on we consider homology with coefficients in a fixed finite field $\F$.

For time series segments $\gamma = (x_k,x_{k+1},\dots, x_\ell)$, this criterion cannot be applied directly, even if $x_\ell=x_k$, since a finite sample carries the discrete topology.
Furthermore, we are not only interested in periodic segments, but also in segments that are almost-periodic. For example, if the underlying deterministic dynamical system is subject to noise, a periodic segment might be perturbed into an almost-periodic segment by the noise.  

Using ideas from topological data analysis, we now derive a meaningful generalization of the criterion $H_1(\zeta)\neq 0$. 
For this, we consider \emph{thickenings} of $\gamma$. In general the \emph{$r$-thickening} of a subset $G$ of a metric space $(M,d)$ is defined by
\begin{equation}
	O_r(G) = \{ x \in M\mid d(x,G) < r \}, \qquad 0 \leq r.
	\label{eq:thickening}
\end{equation}
The thickened set $ O_r(\gamma) $ in $X$ can be interpreted as the set that contains all segments $\tilde\gamma$ that are '$r$-perturbations' of $\gamma$ with respect to the metric $d$.

\begin{exmp}

\begin{figure}[h]
	\centering
	\includegraphics[width=0.5\textwidth]{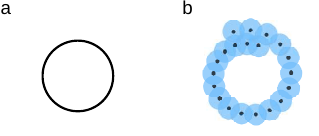}
    {\phantomsubcaption\label{fig:circle}}
	{\phantomsubcaption\label{fig:almost-circle-points-thickening}}
	\caption{a) A periodic orbit is a topological circle.  b) Almost periodic time series segment (black) with thickening (blue). Only the thickening has nontrivial $H_1$.
    \label{fig:cycling-definition-illustrations}
	}
\end{figure}
Let $\gamma$ be the time series segment shown as black dots in \cref{fig:almost-circle-points-thickening}. This segment represents what we want to detect as cycling. 
Clearly, $H_1(\gamma) = 0$. 
However, thickening $\gamma$ with respect to the Euclidean metric (indicated by the blue balls in \cref{fig:almost-circle-points-thickening}), we detect nontrivial homology 
$H_1(O_r(\gamma))\neq 0$ for $r$ in some interval $(r_0,r_1)$.
\qed
\end{exmp}

The example demonstrates how thickening almost-periodic segments 
generalizes the criterion $H_1(\gamma) \neq 0$.
However,
in general, a metric which incorporates only spatial information may not be able to distinguish relevant homological features of a time series from features that arise due to sampling.
This is illustrated in the following example.

\begin{exmp}[Squished periodic orbit]
    \label{ex:squished-periodic-no-tangent-bundle}
    Suppose we have an ellipse-shaped periodic orbit
    $E = \{ x\in \R^2 \mid (x_1/a)^2 + (x_2/b)^2 = 1\}$
    with $ b < a $, as the grey curve in \cref{fig:squished-periodic-segment}.

    Using the metric space $\R^2$ with the metric induced by the Euclidean norm, we have
    \[
        H_1(O_r(E)) = \begin{cases}
            \F, & \text{ if } 0\leq r< b,\\
            0, & \text{ else.}
        \end{cases}
    \]
    Even though we consider a periodic orbit, for $ b $ much smaller than $ a $ only correspondingly small thickening radii $r$ lead to nontrivial $H_1$ since the center of the ellipse quickly closes up with increasing $r$. 
    This poses a problem when passing to a time series $\Gamma$ sampled from the periodic orbit (black points in \cref{fig:squished-periodic-segment}), as it may not be possible to distinguish a relevant homological feature from a feature induced by sampling.
    Many thickened segments $\gamma$ of $\Gamma$ have holes
    and thus fulfill the criterion $H_1(O_r(\gamma)) \neq 0 $.
    For example, we would erroneously detect the segment consisting of the seven leftmost points in the sampling as `cycling'.
    This should be regarded as sampling artifact, since, in this time series, only segments that wrap around the whole ellipse should be detected as `cycling'.
    \qed
\end{exmp}

\begin{figure}[h]
\centering
    \begin{tikzpicture}[scale=1.4]
        \begin{axis}[axis equal,axis lines=none,ymin=-2,ymax=2]
            \addplot[only marks,mark=*,mark size=8.02pt,color=nice_blue,opacity=.7,domain=0:2*pi,samples=17,variable=\t]({cos(deg(t))},{.2*sin(deg(t))});
        \end{axis}
        \begin{axis}[axis equal,axis lines=none,ymin=-2,ymax=2]
            \addplot[color=gray,opacity=1,domain=0:2*pi,samples=51,variable=\t,line width=2pt]({cos(deg(t))},{.2*sin(deg(t))});
            \addplot[only marks,mark=*,mark size=1.5pt,color=black,opacity=.7,domain=0:2*pi,samples=17,variable=\t]({cos(deg(t))},{.2*sin(deg(t))});
        \end{axis}
    \end{tikzpicture}
    \caption{ Periodic segment (black) with thickening (blue) sampled from a squished periodic orbit (gray).}
    \label{fig:squished-periodic-segment}
\end{figure}

In order to be able to detect the relevant homological features, we lift the sampling data to the unit tangent bundle.

The \emph{tangent bundle} $\T(X)$ of $X\subset \R^d$ is the set of all tuples $ (x,v) \in X\times \R^d$, where $x$ is a point in $X$ and $v$ is the velocity of some curve through $x$. 
Any point $ x\in X $ that is not an equilibrium point of $\Phi$ has an associated nonzero tangent vector $ v(x) = \frac{d}{dt} \Phi(t,x) |_{t=0}$. 

For our purpose, it is more appropriate to consider the direction of motion at each non-equilibrium point (and ignore the speed), thus we normalize the second component and pass to the \emph{unit tangent bundle}
\[
	\UTB(X) = X\times S_2^{d-1} = \{ (x,v)\in \T(X)\mid \norm{v}_2=1 \}.
\]
Let $ X_0 $ be the set of equilibria of $\Phi $ and consider the section
\begin{equation}
\label{eq:section}
	\rho\colon X\setminus X_0 \rightarrow \UTB(X\setminus X_0), \quad x \mapsto (x,v(x)/\|v(x)\|_2)
\end{equation}
of the unit tangent bundle induced by the flow. 
For time series data, we recover $v$ through finite differences (which imposes a condition on the sampling density).
We use the (approximated) section $\rho$ in order to lift the time series data $x_0,\ldots,x_N$ from $X$ to $\UTB(X\setminus X_0)\subset \UTB(X)$, i.e., we consider the lifted time series $\rho(\Gamma):=\{\rho(x_0),\ldots,\rho(x_N)\}$ in $\UTB(X)$.   Then, instead of considering distances between points $x,y\in X$ directly, we measure the distance between $\rho(x)$ and $\rho(y)$ in $\UTB(X)$, which incorporates spatial as well as directional information. To this end we introduce a metric in $\UTB(X)$,
\begin{equation}
\label{eq:UT-metric}
    d_C((x,v),(y,w)) = \max\left\{ \norm{x-y}_2, C \norm{v-w}_2 \right\}
\end{equation}
which is parameterized by a weighting factor $C>0$.
The following example demonstrates that by this lifting to the unit tangent bundle, we recover natural dynamical features of segments.

\begin{exmp}[Squished periodic orbit revisited]
    \begin{figure}[h]
    \centering
    \begin{tikzpicture}[scale=1.5]
        \begin{axis}[axis equal,axis lines=none,ymin=-2,ymax=2]
            \addplot[color=gray,opacity=.5,domain=0:2*pi,samples=51,variable=\t,line width=2pt]({cos(deg(t))},{.2*sin(deg(t))});
            \addplot[only marks,mark=*,mark size=1.5pt,color=black,opacity=.7,domain=0:2*pi,samples=17,variable=\t]({cos(deg(t))},{.2*sin(deg(t))});
            \addplot[samples=17, domain=0:2*pi,
            variable=\t,
            quiver={
            u={-sin(deg(t)) / sqrt( (sin(deg(t)))^2 + (0.2*cos(deg(t)))^2) },
            v={.2*cos(deg(t)) / sqrt( (sin(deg(t)))^2 + (0.2*cos(deg(t)))^2)},
            scale arrows=0.2},->,nice_blue,line width=1pt]({cos(deg(t))}, {.2*sin(deg(t))});
        \end{axis}
    \end{tikzpicture}
    \caption{ Periodic segment (black) with unit tangent vectors (blue) sampled from a squished periodic orbit.}
    \label{fig:squished-periodic-segment-2}
    \end{figure}
    The situation in \cref{ex:squished-periodic-no-tangent-bundle} changes if we lift to the unit tangent bundle (see \cref{fig:squished-periodic-segment-2}). 
    Again, we first consider the periodic orbit $E$. Consider the lifts $(x,v)$ and $(y,w)$ of two points that are opposite on the ellipse. 
    Their distance $d_C$ is the maximum of $\norm{x-y}_2$ and $C\norm{v-w}_2 = 2C$. While the first term may be small if $x\approx (0, \pm a)$, the second term is independent of $a$, therefore any thickening radius $r< 2C$ leads to nontrivial $H_1$. 

    This change has important consequences for the sampled segment. Now, the distance between two consecutive points is small compared to the distance between points opposite on the ellipse. 
    A computation shows that suitable $C$ yields
    \[
        H_1(O_r(\gamma)) = \begin{cases}
            \F, & \text{ if } s/2 \leq r< 2C,\\
            0, & \text{ else,}
        \end{cases}
    \]
    where $s$ is the smallest distance such that two consecutive balls intersect. In particular, we obtain a robust criterion for detecting `cycling' segments.
    Furthermore, no segments that wrap around less than half the ellipse are detected.
    \qed
\end{exmp}

The importance of the lift to the unit tangent bundle is further illustrated in the Lorenz data set.

\begin{exmp}[Cycling segments in the Lorenz system]
    Consider segments in the time series in \cref{fig:lorenz-motivating-example} that wrap around the left wing and  can be thickened, using the Euclidean distance, to have nontrivial degree 1 homology.
    Segments that come close to the equilibrium at the center of the left wing have nontrivial $H_1$ only for a small range of thickening parameters. 
    Moreover, one can generate cycling segments close to the equilibrium with an arbitrarily small range of parameters where $H_1$ is nontrivial.

    In contrast, thickenings in the unit tangent bundle take into account the different tangent directions close to the equilibrium.
    This allows robust detection of segments that cycle around the two equilibria at the center of the wings of the Lorenz attractor. 
\end{exmp}

The lift to $\UTB(X)$ and the metric $d_C$ on $\UTB(x)$ for some suitable $C>0$ allow us to make the following definition.

\begin{defn}
    A segment $ \gamma$ is \emph{$r$-cycling} if $H_1(O_r(\rho(\gamma))) \neq 0 $, where $\rho(\gamma)$ is the $r$-thickening of $\rho(\gamma)$ in $\UTB(X)$ with respect to the metric $d_C$.
\end{defn}

\subsection{Classification of Cycling Segments}

While the highlighted segments in \cref{fig:lorenz-motivating-example} are detected to be $r$-cycling for suitable $r$, we have not yet identified the qualitative similarity and dissimilarity of these segments. In the following, we explain how a comparison space can be used to classify cycling segments.

\begin{defn}
    A \emph{comparison space} for a time series $\Gamma\subset X$ is a neighborhood $Y\subset \UTB(X)$ of $\rho(\Gamma)$ where $\rho$ is the section of $\UTB(X)$ defined in \eqref{eq:section}.
\end{defn}
Note that there exists an $r_0=r_0(Y)>0$ such that for any $r<r_0$ and for any segment $ \gamma $ of $ \Gamma $, we have an inclusion map $ i_{\gamma,r}\colon O_r(\rho(\gamma)) \rightarrow Y $. 
By functoriality of homology, the map $i_{\gamma,r}$ induces a linear map
\[
H_1(i_{\gamma,r})\colon H_1(O_r(\rho(\gamma))) \rightarrow H_1(Y)
\]
in homology.
The image of this map is a subspace $\im H_1(i_{\gamma,r}) $ of $ H_1(Y) $.
We compare different cycling segments by their induced subspaces of the vector space $ H_1(Y) $.

\begin{defn}
    The \emph{$r$-cycling space of} a segment $\gamma$ with respect to the comparison space $Y$ is the subspace
    \[
        \Cyc_r(\gamma,Y) = \im H_1(i_{\gamma,r})
    \]
    of $H_1(Y)$. We write $\Cyc_r(\gamma)$ for $\Cyc_r(\gamma,Y)$ in case the comparison space is clear from the context. The \emph{$r$-cycling rank} is the dimension of $\Cyc_r(\gamma)$. We call $H_1(Y)$ the \emph{homological comparison space}.
    \label{def:r-cycling-space}
\end{defn}

In the following, we illustrate the concept of comparison spaces and the cycling space using a time series from a stochastically perturbed double well system. 

\begin{exmp}[Double well system without unit tangent bundle]
To illustrate this idea, consider the time series in Fig.~\ref{fig:dw-cycling-comparison}.
For clarity, we initially disregard the unit tangent bundle (i.e., $C = 0$ in the metric \eqref{eq:UT-metric}) and let the union of the boxes, denoted by $Z$, take the role of the comparison space.
Since $Z$ has two holes, we have $ \dim H_1(Z) = 2$.

\begin{figure}[h]
    \centering
    \includegraphics[width=.7\linewidth]{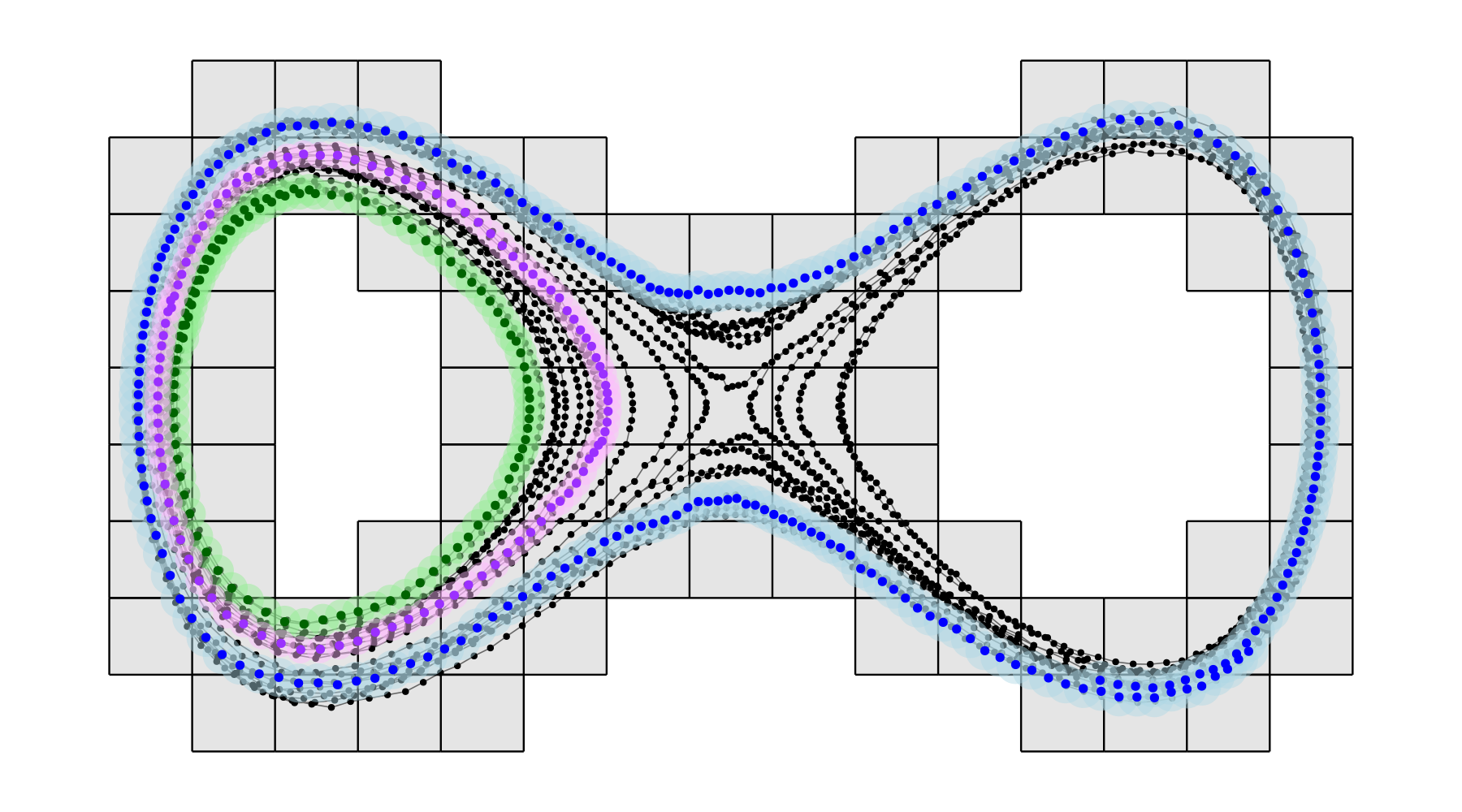}
	\caption{Time series with (simplified) comparison space. The time series (black dots) is covered with boxes, the union of the boxes is the comparison space $Y$. Three cycling segments are highlighted (dark blue, dark green dark purple) together with a thickening (light blue, light green, light purple). The comparison space classifies these into two different classes, one containing the blue, the other containing the green and purple segments. }
	\label{fig:dw-cycling-comparison}
\end{figure}

Three segments of the times series (black) are highlighted: a green segment $\gamma_g$ and a purple one $\gamma_p$ which only loop around the left hole, as well as a blue one $\gamma_b$ which orbits the right and the left hole.

For a suitable range of $r$, the thickenings  of these segments have the topology of a circle and therefore one-dimensional $ H_1 $, where the homology is generated by a curve that loops exactly once around the circle.
By inclusion, these curves also give rise to homology classes in $Z$ and, with slight abuse of notation, we have
\[
	\Cyc_r(\gamma_g, Z) = \Cyc_r(\gamma_p, Z) \neq \Cyc_r(\gamma_b, Z),
\]
Therefore, we detect that the green and purple segments are qualitatively similar while the blue segment is qualitatively different.
\qed
\end{exmp}

\begin{figure}[h]
    \centering
    \includegraphics[width=.7\linewidth]{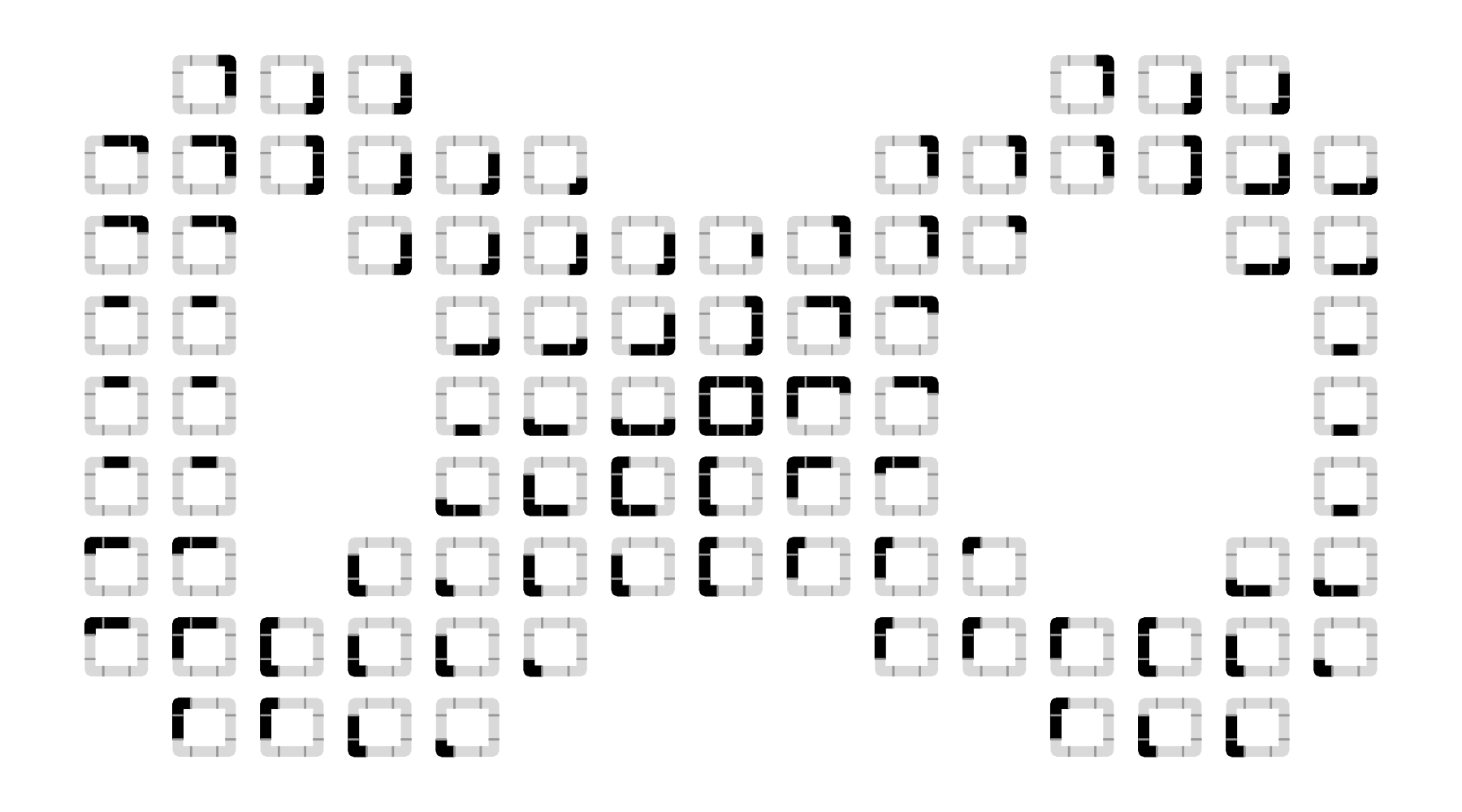}
    \caption{Comparison space in the unit tangent bundle. 
    We use the $S^1_\infty$ unit tangent bundle instead of $S^1_2$ in this figure (and for computations), since $S^1_\infty$ is easier to cover with boxes than $S^1_2$.
    Each cluster of boxes corresponds to the cover of $Q\times S^1_\infty$ where $Q$ is a box in \cref{fig:dw-cycling-comparison}. 
    The comparison space is the union all the black boxes.
    }
    \label{fig:dw-sb-box-space}
\end{figure}

\begin{exmp}[Double well system with unit tangent bundle]
    We now consider the same time series in \cref{fig:dw-cycling-comparison}, but
    lifting the time series to the unit tangent bundle before constructing a cubical cover.     
    In the tangent space, we then consider boxes centered on a grid and take these to cover the lifted data. In our example, we consider boxes of the form $\overline{B_\infty(p,r/2)}\times \overline{B_\infty(q,1/2)}$, where $p\in r\Z^d$ and $q\in \Z^d$.
    
    The resulting comparison space $Y$ in the tangent bundle is visualized in \cref{fig:dw-sb-box-space}.
    Every black box in this figure is to be understood as the product of a box from \cref{fig:dw-cycling-comparison} with a box in the tangent component.
    As we can see, over all but one of the boxes from \cref{fig:dw-cycling-comparison}, a contractible segment of the sphere is covered. 
    Only in the center of the picture, the entire sphere is covered with boxes, since the center box in \cref{fig:dw-cycling-comparison} contains points moving in all directions.
    This additional ``hole'' leads to $\dim H_1(Y) = 3$.
    The double well system (without perturbation) has a hyperbolic equilibrium at the origin, which is precisely where this extra `hole' in the comparison spaces $Y$ is located.
    The appearance of this extra `hole' is therefore expected, since close to a hyperbolic equilibrium, all tangent directions occur.
    Having explained the difference between $Y$ and $Z$, we now use $Y$ to compare the three highlighted segments. Writing $\gamma_g$, $\gamma_p$ and $\gamma_b$ as in the previous example, we have three 1-dimensional cycling spaces with
    \[
	   \Cyc_r(\gamma_g,Y) = \Cyc_r(\gamma_p,Y) \neq \Cyc_r(\gamma_b,Y).
    \]
    As expected, the segments are classified in accordance with our intuition.
    \qed
\end{exmp}

In the previous example, we only considered one-dimensional cycling spaces.
In general, the dimension of a cycling space may be as large as the dimension of the homological comparison space.

\begin{exmp}[Higher-dimensional cycling space]
\label{ex:higher-rank-segments}
\begin{figure}[H]
    \centering
    \begin{subfigure}[b]{0.49\textwidth}
        \centering
        \includegraphics[width=\textwidth]{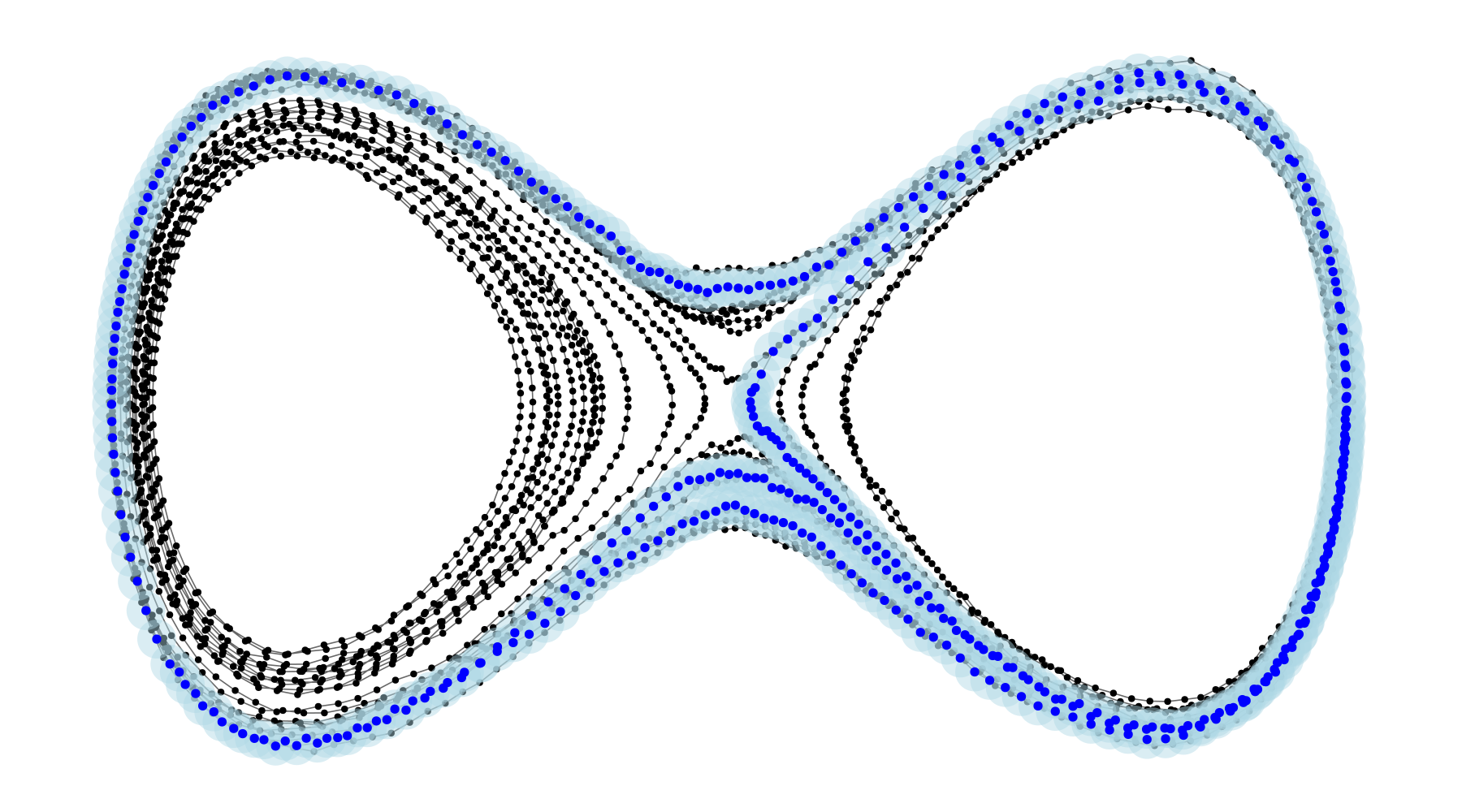}
        \caption{}
        \label{fig:dw-rank-2-segment}
    \end{subfigure}
    \begin{subfigure}[b]{0.49\textwidth}
        \centering
        \includegraphics[width=\textwidth]{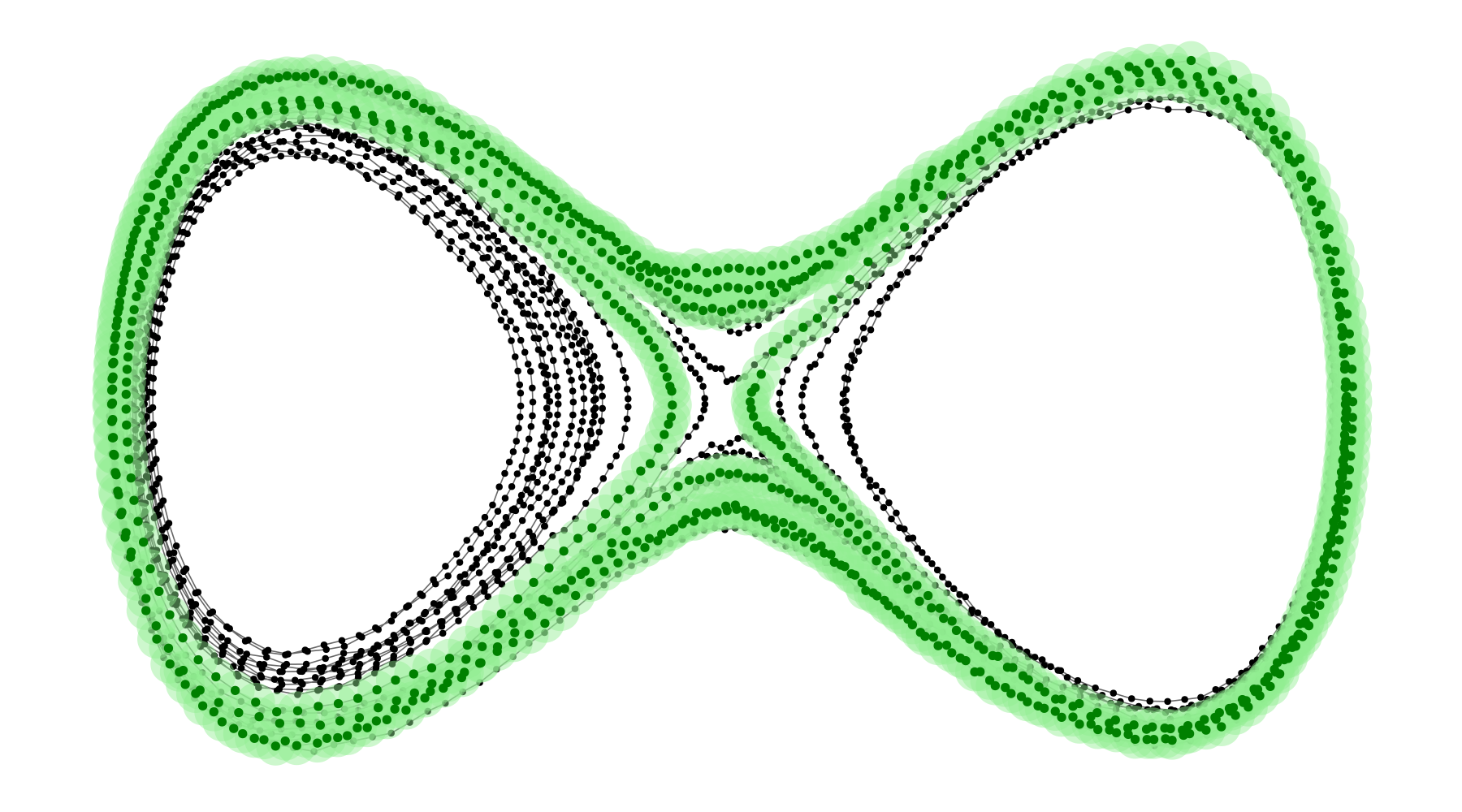}
        \caption{}
        \label{fig:dw-rank-3-segment}
    \end{subfigure}
    \caption{Segments with higher-dimensional cycling spaces.}
    \label{fig:dw-higher-rank-segments}
\end{figure}
    We consider the segments shown in blue and green in \cref{fig:dw-higher-rank-segments},  which we denote by $\gamma_{b2}$ and $\gamma_{g2}$, respectively. 
    
    For exposition, we again start in the simplified setting where the union of boxes $Z$ in \cref{fig:dw-cycling-comparison} takes the role of the comparison space.
    In this case, the cycling spaces of both segments in \cref{fig:dw-higher-rank-segments} are two-dimensional and we have
    \[
        \Cyc_r(\gamma_{b2},Z) = \Cyc_r(\gamma_{g2},Z) = H_1(Z)
    \]  
    Note that we do not detect the qualitative difference of the segments: $\gamma_{b2}$ does not contain a subsegment which wraps only around the left hole in the dataset, while $\gamma_{g2}$ does. This is an instance of the general fact that a vector spaces does not have a canonical decomposition into lower-dimensional subspaces. 
    
    Using the comparison space $Y$ in the unit tangent bundle, the cycling space $\Cyc_r(\gamma_{b2},Y) $ has dimension 2 while $ \dim \Cyc_r(\gamma_{g2},Y) = 3$. Indeed, the cycling space of $\gamma_{b2}$ does not contain the cycling space of segment $\gamma_{p}$ from \cref{fig:dw-cycling-comparison} as a subspace since no curve in the blue set wraps around the left hole but not the center hole.
    
    While the cycling spaces alone do not allow us to conclude that $\gamma_{g2}$ contains a subsegment with cycling space $\Cyc_r(\gamma_p,Y)$, we can conclude from $\Cyc_r(\gamma_p,Y) \not\subset \Cyc_r(\gamma_{2b},Y)$ that 
    $\gamma_{b2}$ does not contain a subsegment with that cycling space.
    \qed
\end{exmp}

\begin{rem}[Comparison Space for the Lorenz time series]
    We return to the Lorenz time series from \cref{fig:lorenz-motivating-example}.    
    A comparison space to detect and distinguish the cycling segments on the wings should have a hole at the center of each wing.
    Covering the time series directly, i.e., without lifting to the unit tangent bundle, would require impractically small boxes that result in a box cover with a very large amount of boxes. 
    Such a cover would furthermore form a hole on the outside of the right wing as an artifact of the time series not reaching this part of the attractor. 

    In contrast, the lifted data set can be covered with rather coarse boxes without losing the holes at the centers of the wings.  
    This is because a closed curve can have arbitrarily small length in space, while its length in the unit tangent bundle is bounded from below by $2\pi C$ (as a consequence of Fenchel's theorem \cite{Fe29}, which states that the total curvature of any closed curve is at least $2\pi$).
    The unit tangent bundle is therefore crucial in order to construct a comparison space that allows us to distinguish the different classes of cycling segments.
\end{rem}

\subsection{The Cycling Signature}
A priori, it is not clear how to appropriately choose the thickening radius $r$ in the construction of the cycling space.
Clearly, the cycling spaces of the segments under consideration will depend on $r$, and different choices of $r$ may be relevant depending on the scales in the system. We adopt the viewpoint of persistent homology \cite{Ghrist.2007,edelsbrunner2008,carlsson2009,MR2121296} 
and argue that we should be interested in 
the entire range of possible thickening radii. To formalize this, we recall the following notion.

\begin{defn}
    Let $I\subset\R$. A \emph{filtered vector space} $V$ is a collection of vector spaces $(V_r)_{r\in I}$ such that $V_r\subset V_s$ for $r\leq s$, together with inclusion maps $V_{r,s}\colon V_r \rightarrow V_s$.
\end{defn}

As we show next, the family of $r$-cycling spaces $\Cyc_r(\gamma,Y)$, $r\in [0,r_0(Y))$, forms a filtered vector space. For this, let $\gamma$ be a segment of $\Gamma$ and $Y$ a comparison space for $\Gamma$. If $r\leq s < r_0(Y)$, we have the inclusions $O_r(\rho(\gamma))\subset O_s(\rho(\gamma))\subset Y$.
Thus, the inclusion map $i_{\gamma,r}\colon O_r(\rho(\gamma)) \rightarrow Y$ factorizes as $i_{\gamma,r} = i_{\gamma,s} \circ O_{r,s}(\rho(\gamma))$, where $O_{r,s}(\rho(\gamma))$ is the inclusion map $ O_r(\rho(\gamma)) \rightarrow O_s(\rho(\gamma))$. We thus have
\[
    \Cyc_r(\gamma,Y) = \im H_1(i_{\gamma,r}) = \im H_1(i_{\gamma,s} \circ O_{r,s}(\rho(\gamma))) \subset \im H_1(i_{\gamma,s}) = \Cyc_s(\gamma,Y).
\]
Therefore, the cycling spaces form a filtered vector space, and we can give the following definition.

\begin{defn}
    \label{def:cycling-signature}
    The \emph{cycling signature} of a segment $ \gamma $ of $\Gamma$ in $Y$ is the filtered vector space 
    \[ 
    \Cyc(\gamma,Y) = (\Cyc_r(\gamma,Y))_{r\in I}, 
    \]
    where $I = [0,r_0(Y))$  is the \emph{interval of admissible radii}.
    The \emph{cycling rank function} of $\gamma$ is the function 
    \[
    r \mapsto \dim \Cyc_r(\gamma,Y).
    \]
\end{defn}

\begin{exmp}[Rank and thickening radius]
\begin{figure}[H]
    \centering
    \begin{subfigure}[b]{0.49\linewidth}
        \centering
        \includegraphics[width=\textwidth]{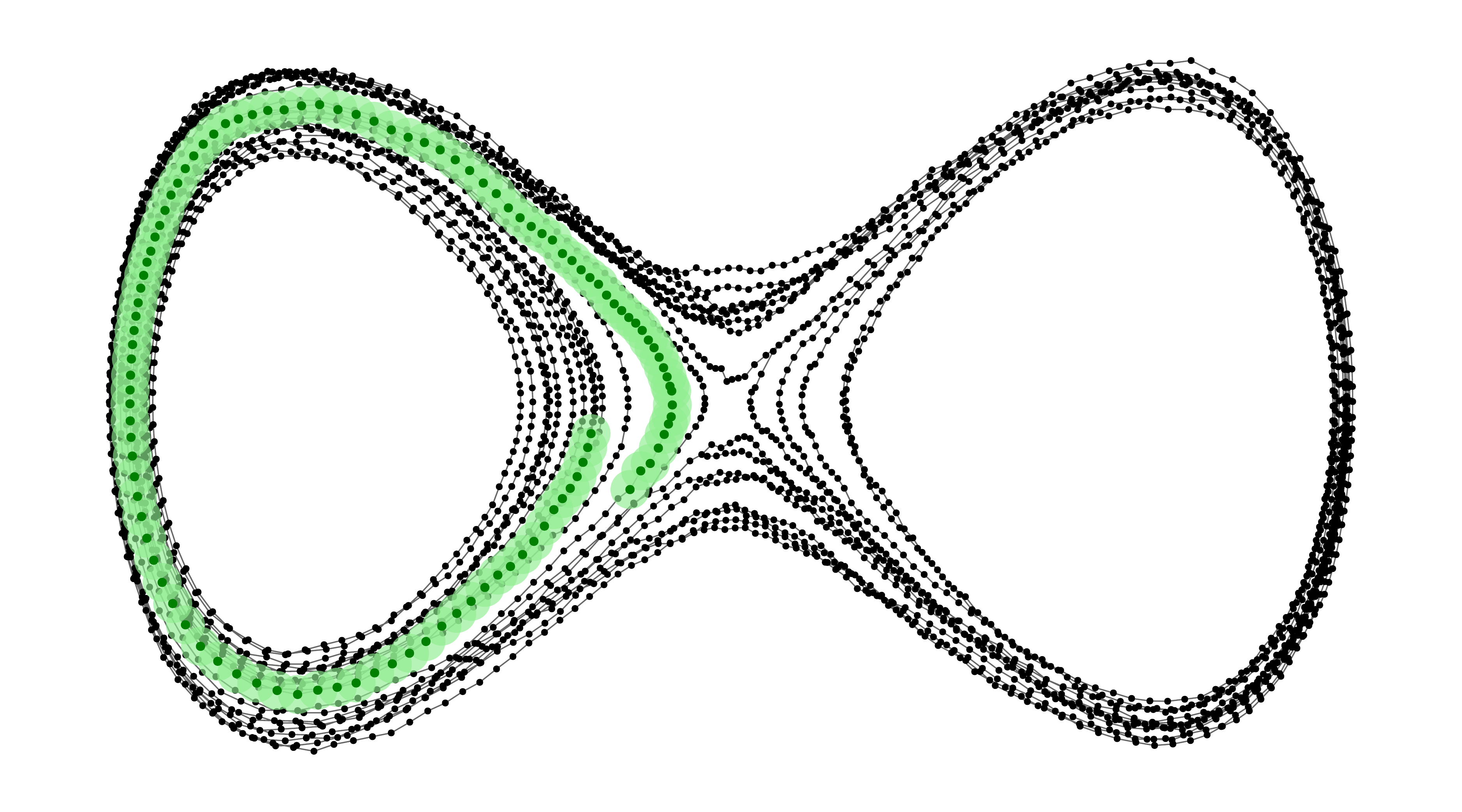}
        \caption{}
        \label{fig:dw-rank-change-0}
    \end{subfigure}
    \begin{subfigure}[b]{0.49\linewidth}
        \centering
        \includegraphics[width=\textwidth]{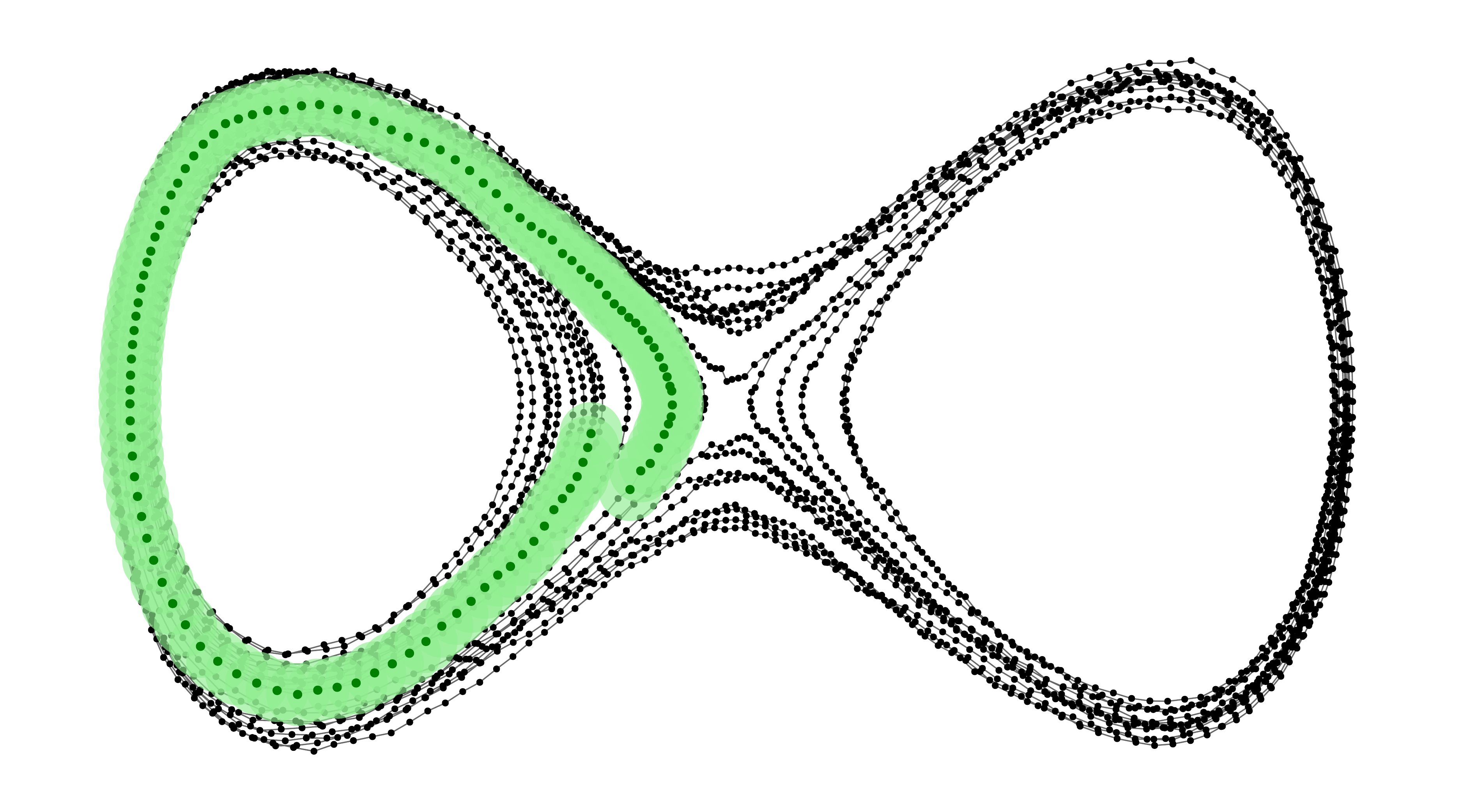}
        \caption{}
        \label{fig:dw-rank-change-1}
    \end{subfigure}
    \caption{Segment with different thickenings. The light green set in \textbf{a} has rank 0, and the one in \textbf{b} has rank 1.}
    \label{fig:dw-rank-change}
\end{figure}

    Consider the segment in \cref{fig:dw-rank-change} and the comparison space $Y$ from \cref{ex:higher-rank-segments}. From the figure, it is clear that small thickenings lead to rank 0, while larger thickenings lead to rank 1.
    This is captured in the cycling rank: there is a threshold $\bar{r}$ such that
    \[
        \dim \Cyc_r(\gamma, Y) = \begin{cases}
            0,& \text{ if } r < \bar{r},\\
            1,& \text{ if } \bar{r} \leq r < r_0(Y).
        \end{cases}
    \]
    \qed
\end{exmp}

\section{Statistics of Cycling: Computational Experiments}
\label{sec:application}

The notions on cycling introduced in \cref{sec:math} apply to a single segment $\gamma=(x_k,\ldots,x_\ell)$ of a given time series $\Gamma\subset\R^d$.
In order to obtain an overview on the cycling properties of different segments from a (long) time series $\Gamma$, we carry out the following sampling approach: 
We choose a (finite) set of time spans $\mathcal{T}\subset [0,\infty)$, a sample size $N\in\N$ and randomly sample $N$ segments of time span $\tau$ for each $\tau\in \mathcal{T}$. We then compute the cycling signature for these $N\cdot |\mathcal{T}|$ segments as described in \cref{sec:computation} and analyze the resulting collection of signatures.

In the following, we perform this experiment on three examples (for details see \cref{sec:data-generation}): the Lorenz system with the classical parameters \cite{Lo-63}, a perturbed Hamiltonian system derived from a scalar potential with two local minima (the \emph{double well system}) and the Dadras system \cite{Dadras.2012}, a four-dimensional ordinary differential equation with a hyperchaotic attractor. 

\begin{figure}[h]
	\centering
	\includegraphics[width=.32\linewidth]{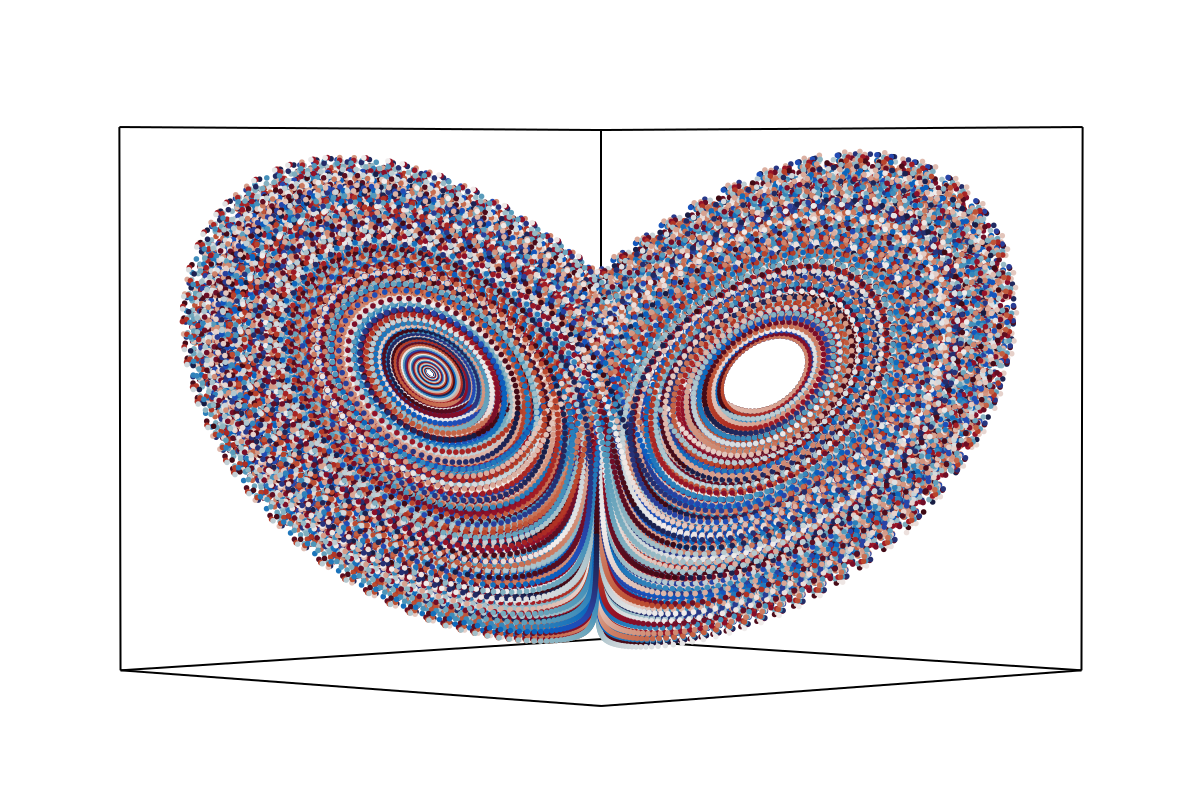}
	\includegraphics[width=.30\linewidth]{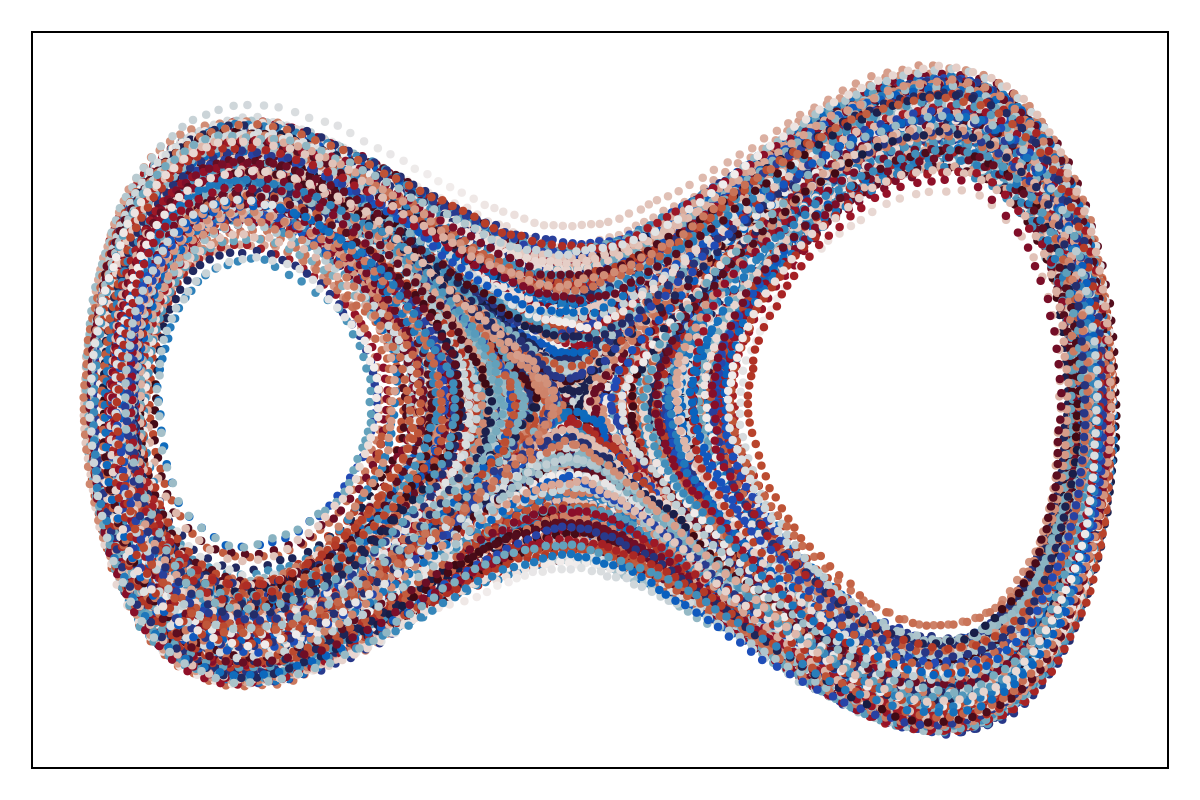}
	\includegraphics[width=.32\linewidth]{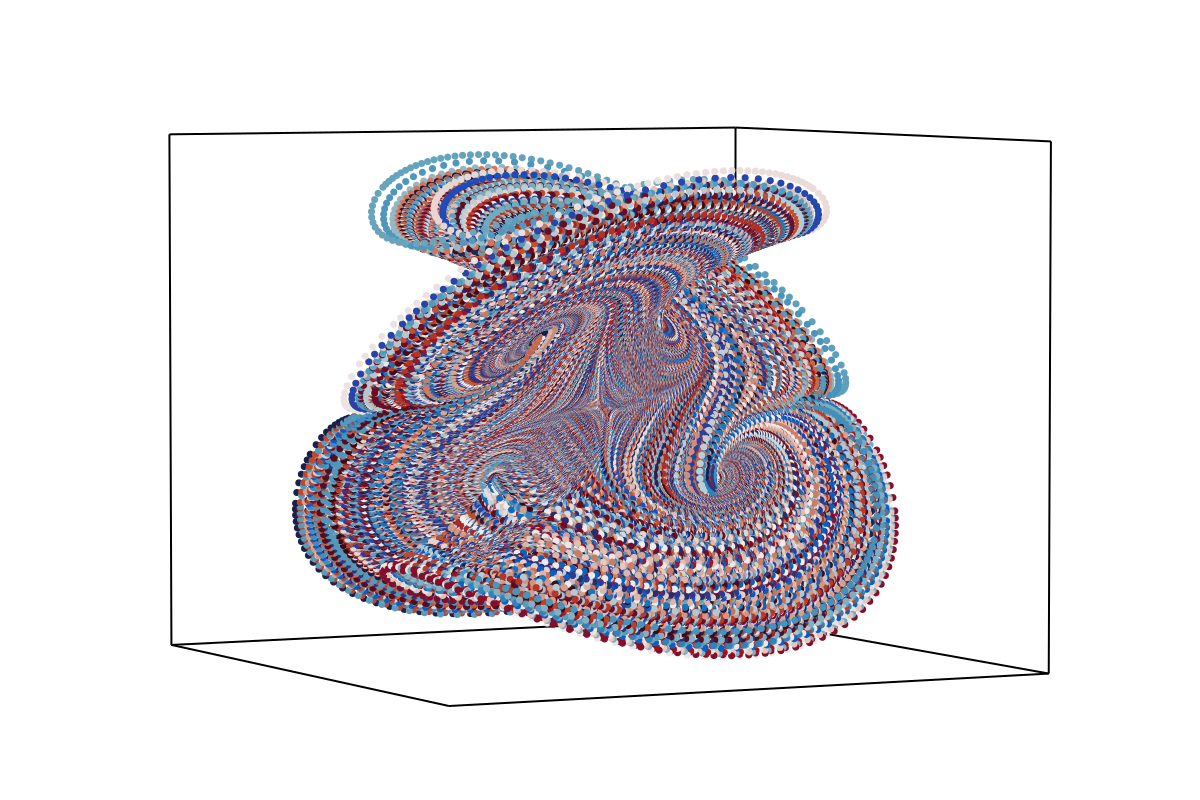}
	\caption{Time series for the Lorenz system (left), a stochastically perturbed Hamiltonian system with a double well potential (center) and the Dadras system projected to the first three coordinates (right). The color of the $x_i$ changes gradually with time $t_i$ to visualize time evolution.}  
	\label{fig:time-series-lorenz-double-well-dadras}
\end{figure}

For each system, we generate a time series $\Gamma$ by numerical integration, cf.\ \cref{fig:time-series-lorenz-double-well-dadras}, and compute a comparison space. This yields a 2d, 3d and 5d homological comparison space for the Lorenz, double well, and Dadras systems, respectively.
We then sample $N=1000$ segments of time spans $\tau \in \{ 10,20,\dots, 1000 \}$ and compute their cycling signatures.

It is not necessary to know the precise construction of the comparison space (which will be described in detail in \cref{subsec:cubical-comparison-space}) to understand the content of this section. However, we do mention that we take boxes of a certain size $(p,k)$, where for Lorenz $(p,k) = (8,2)$, for double well $(p,k) = (0.2,2)$ and for Dadras $(p,k) = (4,3)$. With $C = p/k$, our heuristics (see \cref{subsec:heuristics-C-r0}) suggest $r_0(Y) = p$ as maximal filtration value for the cycling signatures. The different choices of box size are due to the different scales present in the three example systems.

\subsection{Analysis of Cycling Spaces for a Fixed Thickening Radius}

As mentioned in \cref{sec:math}, the cycling space of a segment depends on the thickening radius. 
In this section, for each experiment, we fix a thickening radius $r$ ($r=6$ for Lorenz, $r=1.8$ for double well and $r=3.8$ for Dadras) and analyze the collection of cycling spaces resulting from evaluating the cycling signatures of the sampled segments at that radius, i.e.
the collection of 10000 vector spaces $\Cyc_{r}(\gamma,Y)$. 

Fixing a thickening radius corresponds to fixing a scale of interest. In \cref{subsec:analysis-with-radius}, we will explain our choice of thickening radii and how the results depend on them.  

\subsection*{Cycling rank}

\begin{figure}[H]
    \centering
	\includegraphics[width=0.7\textwidth]{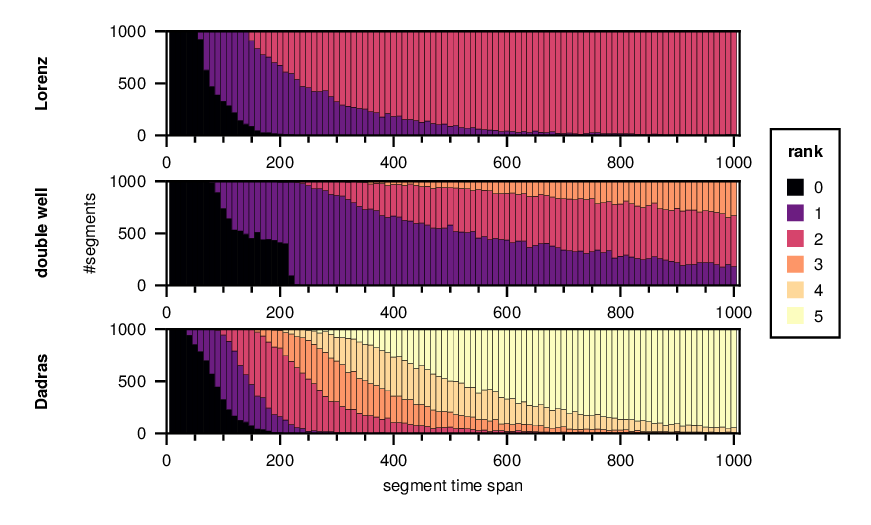}
	\caption{Distribution of cycling ranks for segments of different time spans. Each plot consists of stacked barplots (one for each value of $\tau$) showing the distribution of cycling ranks for the sampled segments. In all three plots, rank zero dominates for short time spans while long segments eventually attain the maximal possible cycling rank of the system.
    }
	\label{fig:cycling-ranks-combined}
\end{figure}

The distribution of cycling ranks as a function of the time span $\tau$ is shown in \cref{fig:cycling-ranks-combined}.  In all three plots, any sufficiently long segment has positive cycling rank, 
i.e., any sufficiently long segment of the time series is cycling.
Furthermore, the time span at which the number of rank 0 segments begins to decrease suggests the minimal time span of a cycling segment (Lorenz $\sim 60$; double well $\sim 80$; Dadras $\sim 40$) and the time span at which there are no rank 0 segments left indicates how long segments can be at most before they exhibit cycling (Lorenz $\sim300$; double well $\sim230$; Dadras $\sim200$).

\subsection*{Rank 1 signatures}

\begin{figure}[H]
    \centering
	\includegraphics[width=.7\linewidth]{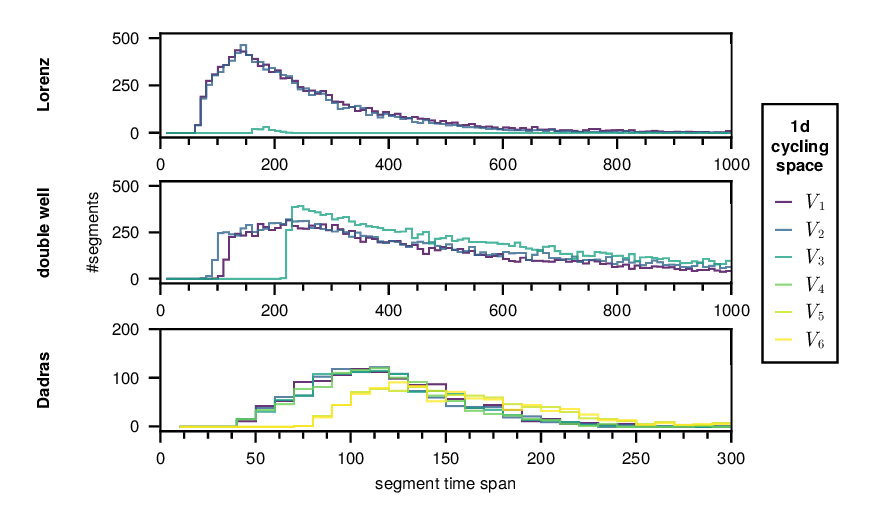}
	\caption{Distribution of rank 1 segments for different time spans. For the double well and Dadras systems, some infrequent signatures are omitted.} 
	\label{fig:1d-subspaces-combined}
\end{figure}

Segments of rank 1 are of particular interest as they exhibit a single type of cycling, which is uniquely identified by its cycling space. 
\cref{fig:1d-subspaces-combined} shows the distributions of rank 1 segments in the example time series.  We observe three different rank 1 cycling spaces in the Lorenz system, three in the double well system, and six in the Dadras system.
Note that this observation does not conflict with the fact that the comparison spaces for the Lorenz and the Dadras system have dimensions 2 and 5, respectively, as the rank 1 signatures may be linearly dependent. 

The smallest time span $\tau_{min}(V)$ for which a particular one-dimensional cycling space $V$ is sampled is an upper bound on the minimal time span of any cycling segment in this space. 
Both the Lorenz and the double well system possess two cycling spaces $V_1, V_2$ with small $\tau_{min}$ compared to a third cycling space $V_3$ (Lorenz: $\tau_{min}(V_1) \approx \tau_{min}(V_2) \approx 70$, $\tau_{min}(V_3) \approx 160$; double well: $\tau_{min}(V_1) \approx 70,\, \tau_{min}(V_2) \approx 90$, $\tau_{min}(V_3) \approx 210$), whereas in Dadras, four spaces $V_1,V_2,V_3,V_4$ have small $\tau_{min}$ compared to $V_5,V_6$.

\begin{figure}[H]
	\centering
	\includegraphics[width=\linewidth]{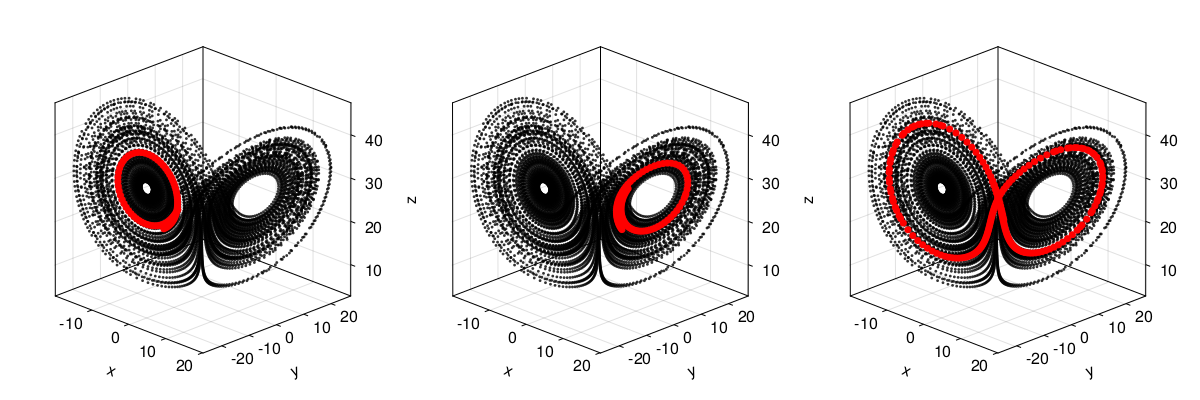}
	\caption{Rank 1 segments in the Lorenz system. The full time series $\Gamma$ is shown in black. The red points in each plot highlight a segment with cycling space $V_1$, $V_2$ and $V_3$, respectively.} 
	\label{fig:lorenz-rank-1-segments}
\end{figure}

In the Lorenz system this is expected, since the third type of cycling, the longer one, is the "sum" of the two shorter ones on the individual wings and traces out both wings as shown in \cref{fig:lorenz-rank-1-segments}.

\subsection*{Higher rank segments.}

More complicated cycling behavior is captured by segments with $r$-cycling rank greater than one.
As mentioned in \cref{ex:higher-rank-segments}, the cycling space itself does not reveal which rank 1 cycling spaces give rise to it, since vector spaces do not canonically decompose into one-dimensional subspaces. 
However, higher dimensional cycling spaces can still be used to relate lower dimensional cycling spaces with each other via subspace inclusion.

\begin{figure}[H]
    \centering
	\includegraphics[width=.7\linewidth]{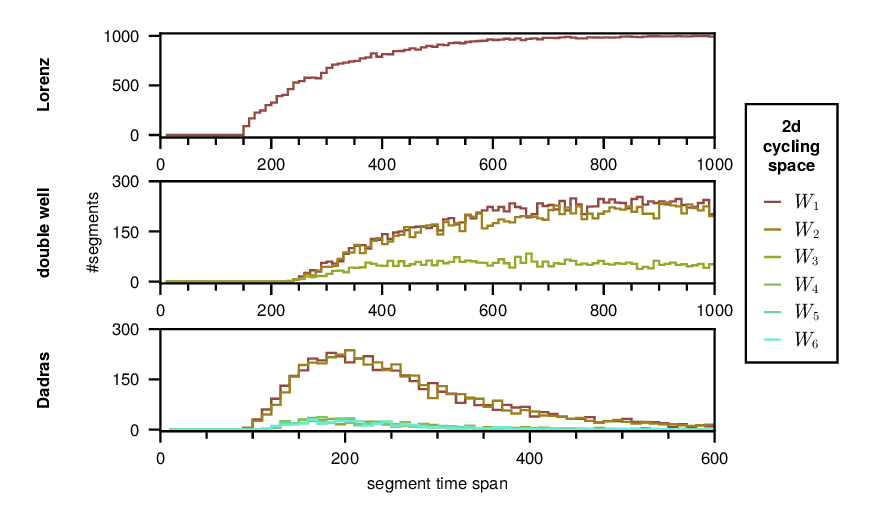}
	\caption{Distribution of rank 2 signatures in the example systems.}
	\label{fig:2d-subspaces-combined}
\end{figure}

Of particular interest is the case of rank 2 segments.  In \cref{fig:2d-subspaces-combined} we show the distribution of rank 2 segments as a function of the segment time span. 
A segment which consecutively visits exactly two different 1d cycling spaces has a 2d cycling space which is the direct sum of the two 1d spaces.
Conversely, there cannot be direct transitions between two specific 1d cycling spaces if there is no segment with a 2d cycling space that contains these two 1d spaces.
This allows us to rule out direct transitions between certain cycling spaces. 

\begin{figure}[H]
    \centering
	\includegraphics{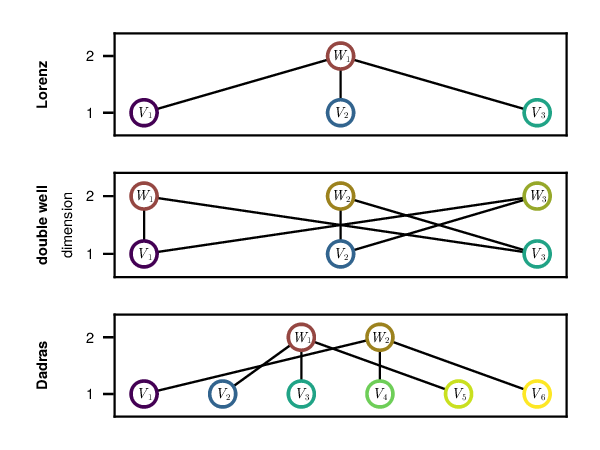}
	\caption{Inclusion graphs. In each plot, the nodes in the bottom and top rows correspond to frequent rank 1 and 2 signatures, respectively. Node colors correspond to the ones in \cref{fig:1d-subspaces-combined} and \cref{fig:2d-subspaces-combined}. An edge indicates that the rank 1 signature is a subspace of the rank 2 signature.}
	\label{fig:s21-inclusion}
\end{figure}

In the Lorenz system, the homological comparison  space is two-dimensional, therefore all rank 2 segments have the same cycling space, and all three of the frequent 1d cycling spaces are subspaces of this space (see \cref{fig:s21-inclusion}). 
Algebraically, the longer cycling space $V_3$ is therefore the span of the two shorter ones $V_1$ and $V_2$, cf.\ \cref{fig:lorenz-rank-1-segments}. 

In the other two systems, rank 2 segments provide deeper insight into the transitions between cycling spaces:
We cannot rule out any transition between the three prominent cycling spaces in the double well system. 
To see this, note that the statistics of the 2d cycling spaces, depicted in \cref{fig:2d-subspaces-combined}, reveal three prominent 2d cycling spaces. 
It can be seen in \cref{fig:s21-inclusion} that any pair of 1d cycling spaces is contained in a 2d cycling space (for example $V_1$ and $V_2$ are contained in $W_3$).

For transitions that occur we obtain quantitative statements:
As can be seen from \cref{fig:2d-subspaces-combined}, the spaces $W_1$ and $W_2$ appear much more frequently than $W_3$.
The inclusion graph in \cref{fig:s21-inclusion} (center) reveals that $W_1$ and $W_2$ contain $V_3$, the space with large $\tau_{min}$, and one of the others $V_1$ and $V_2$, whereas $W_3$ contains both of the spaces with small $\tau_{min}$. 
An interpretation of this is that transitions between $V_1$ and $V_2$ are less likely than transitions involving $V_3$.

In the Dadras system, we detect two prominent 2d cycling spaces which organize the six 1d cycling spaces into two clusters. 
These are defined by $W_1$ and $W_2$, and consist of $V_1,V_4,V_6 \subset W_1$ and $V_2,V_3,V_5 \subset W_2$, i.e., each contains two subspaces with smaller, and one with larger $\tau_{min}$. 
Note that $W_1$ and $W_2$ are disjoint, therefore direct transitions between cycling spaces in the same cluster are frequent while direct transitions between cycling spaces from different clusters are rare.  

\subsection*{Distinguishing the three systems via cycling}

Visually, there are similarities between the time series of the Lorenz and the double well system (see \cref{fig:time-series-lorenz-double-well-dadras});
they both have two centers (`holes') around which trajectories cycle.
The dimensions of the comparison cycling spaces differ (2 for Lorenz, 3 for double well), but in and of itself this does not necessarily imply significant differences in the dynamics. 
However, we can quantify differences using the statistics of cycling signatures.

\cref{fig:1d-subspaces-combined} shows that the Lorenz system has two prominent 1d cycling spaces with similar statistics, suggesting both cycling motions occur with similar frequency.
The double well system also has two prominent 1d cycling spaces, albeit with different $\tau_{min}$.
In both systems, there is a third 1d cycling space with $\tau_{min}$ approximately twice as large as the others.
In the Lorenz system, this third space is observed infrequently, whereas in the double well system, it dominates in frequency once it appears at a time scale of $\tau_{min} \approx 210$. 

As indicated in \cref{fig:cycling-ranks-combined}, the Dadras system possesses six cycling motions while the other two systems only have three.
In addition, the statistics of 1d cycling spaces in \cref{fig:1d-subspaces-combined} show that both Lorenz and double well have clear thresholds where 1d cycling spaces appear, while in Dadras no such threshold can be determined. 

\subsection*{Cycling Signatures in the Dadras system.}

The Dadras system is a 4d system of ordinary differential equations, which appears to exhibit a (hyper)chaotic attractor \cite{Dadras.2012}. 
As the system is four-dimensional, visualization is limited to projections that necessarily omit some spatial information.
Cycling signatures allow us to analyze the dynamics in the Dadras system without visualizations.

As stated above, \cref{fig:1d-subspaces-combined} indicates six 1d cycling spaces.
Our analysis of rank two segments (in particular \cref{fig:s21-inclusion}) reveals that these 1d cycling spaces are organized in two Lorenz-like clusters, 
each of which comprises two spaces with smaller $\tau_{min}$ than the third which they span.

Together, the $V_i$, $i=1,\dots,6$, span a four-dimensional subspace of the homological comparison space.
Similar to the double well system, there is a hyperbolic equilibrium at the origin that generates a nontrivial homology class in the homological comparison space. This generator, together with the four 1d cycling spaces $V_1,V_2,V_3,V_4$, generates the 5 dimensional homological comparison space.  

\subsection{Analysis of Cycling Signatures}
\label{subsec:analysis-with-radius}

In the previous section, we analyzed the collection of cycling spaces that arises by evaluating the computed cycling signatures at a fixed evaluation radius. In the following we discuss visualizations which take the radius into account and discuss our choice of evaluation radius.

\begin{figure}[H]
    \centering
    \includegraphics[width=.9\linewidth]{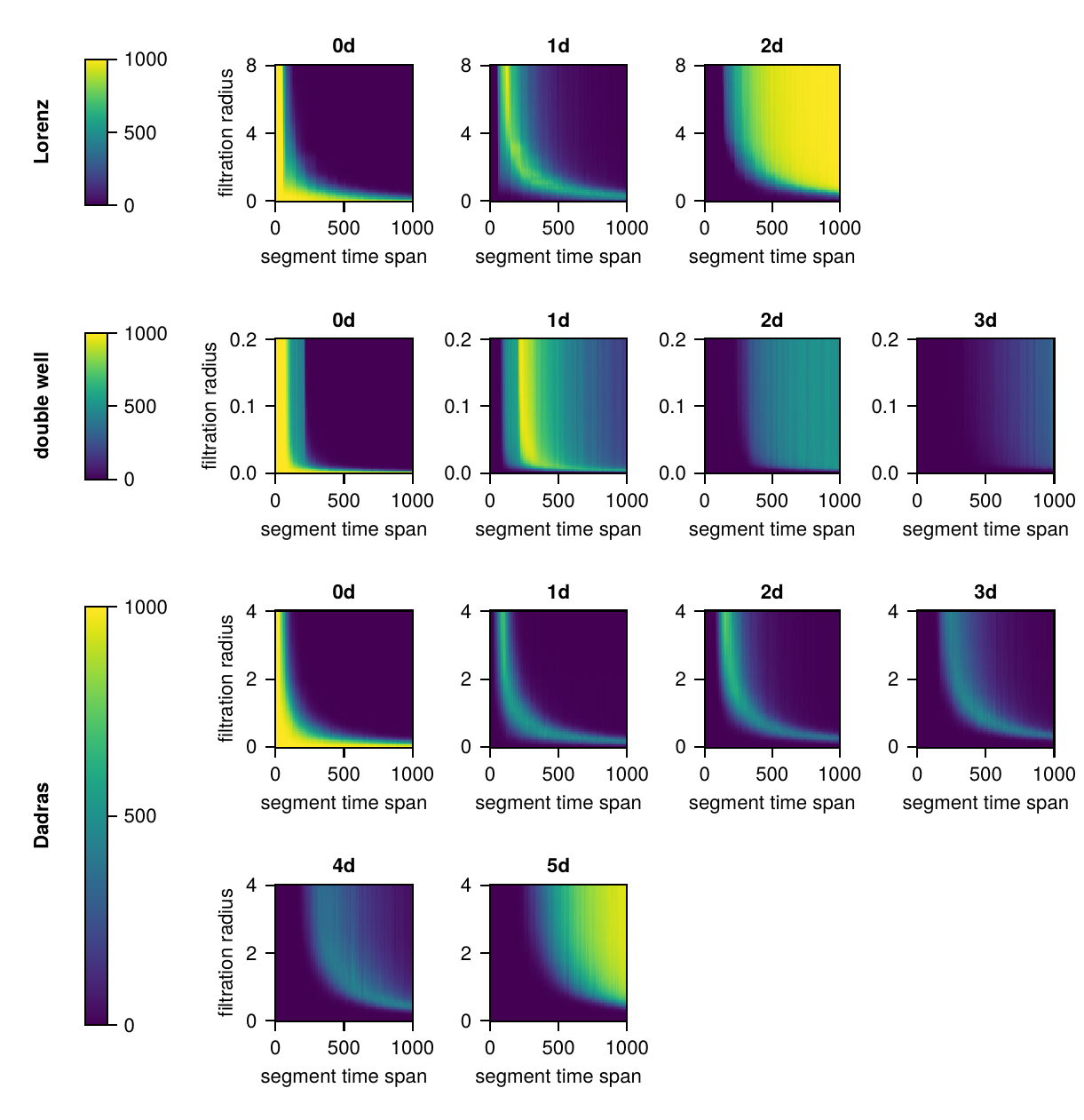}
    \caption{Distribution of cycling ranks for different segment lengths and evaluation radii. Each column corresponds to a dimension, each row to one of the example systems, i.e. the top left heat map corresponds to rank 0 segments in Lorenz.}
    \label{fig:rank-heatmaps}
\end{figure}

\begin{figure}[H]
    \centering
    \includegraphics[width=.9\linewidth]{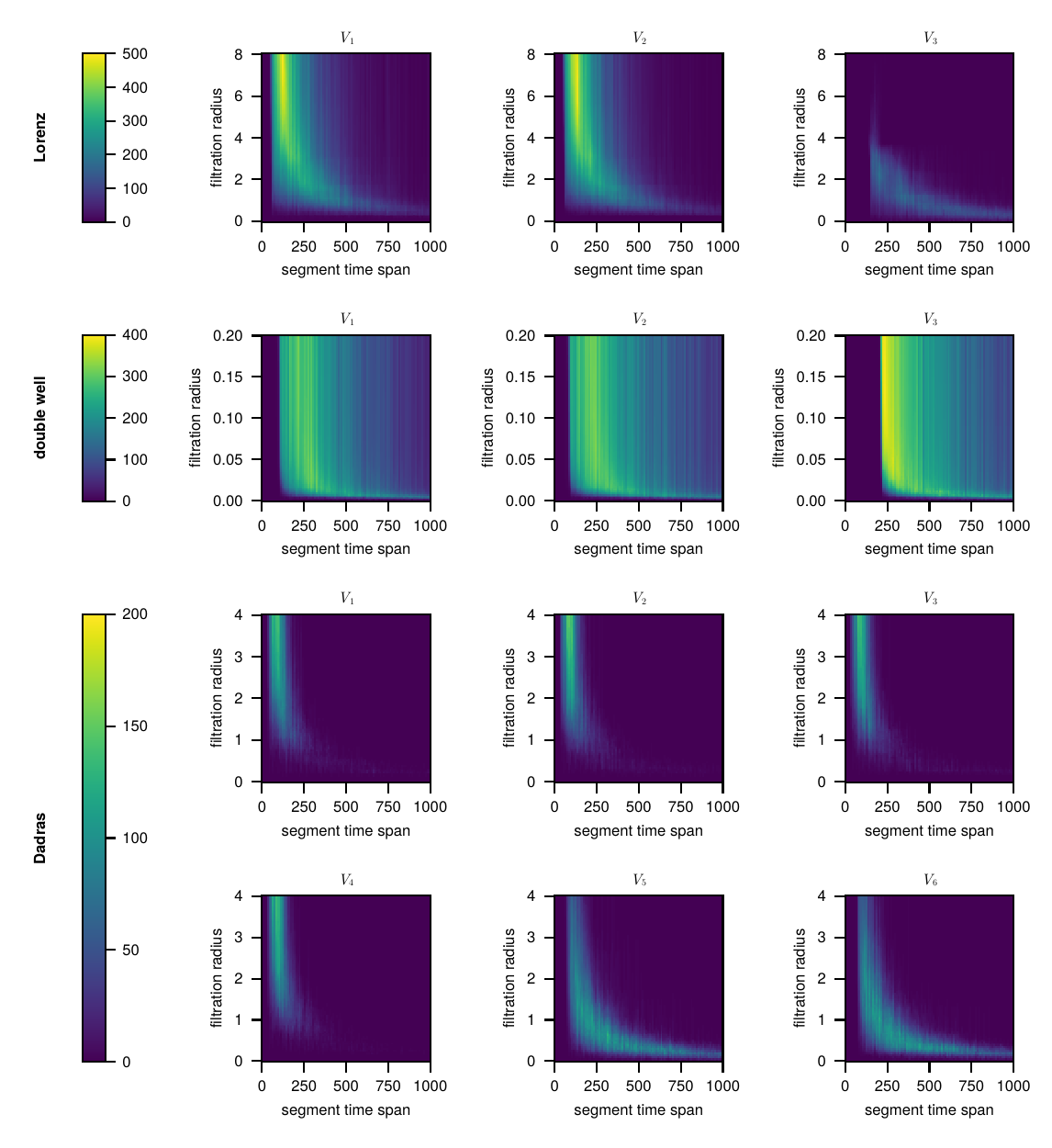}
    \caption{Distribution of the most frequent 1d cycling spaces for different segment lengths and filtration radii. }
    \label{fig:sig1-heatmaps}
\end{figure}

The collection of cycling ranks is visualized in \cref{fig:rank-heatmaps}, the distribution of 1d cycling spaces in \cref{fig:sig1-heatmaps}. The figures in the previous section arise from these heat maps by considering a fixed horizontal slice at a fixed thickening radius.

Qualitatively, not much changes in the top half of the heat maps, which means that there is a large interval of thickening radii that produce similar statistics. This justifies our choice of a fixed radius for each system for the discussion in the previous section. 


\section{Computation of Cycling Signatures}
\label{sec:computation}
In this section, we describe an algorithm to compute the cycling signature of a segment $\gamma$ of a time series $\Gamma$. 
The cycling signature is a filtered vector space over a fixed finite field $\F$.
We thus compute a basis in the following sense.

\begin{defn}
	Let $V=(V_r)_{r\in I}$ be a filtered vector space. A \emph{basis} $B$ of $ V$ is a collection of subsets $ B_r \subset V_r $, $r\in I$, such that
	\begin{itemize}
		\item $B_r \subset B_s$ for $r\leq s$, and
		\item $ B_r $ is a basis for $ V_r $ for all $r\in I$.
	\end{itemize}
\end{defn}

To compute such a basis, we first express the cycling signature as a morphism of persistence modules (\cref{subsec:ph-and-cycling-signature}), which is used for computation. A more detailed overview on the algorithm is then given at the end of \cref{subsec:ph-and-cycling-signature}.

\subsection{Persistent Homology and the Cycling Signature}
\label{subsec:ph-and-cycling-signature}
We start by recalling some basic notions from persistent homology. For a more detailed exposition, see \cite{edelsbrunner2008}.

\begin{defn}
    A \emph{filtered topological space} $X$ indexed by $I\subset \R$ is a collection of topological spaces $(X_r)_{r\in I}$ such that $X_r\subset X_s$ for $r\leq s$. The corresponding inclusion maps are denoted by $X_{r,s}:X_r\to X_s$, $r,s\in I$.
\end{defn}

This definition is analogous to the definition of a filtered vector space. Importantly, metric thickenings $O(G) = (O_r(G))_{r\in I}$ yield a filtered topological space. In the following, we will also use the notions of filtered simplicial complexes and filtered chain complexes without explicitly defining them.

\begin{defn}
    A morphism $\phi\colon X\rightarrow Y$ between filtered topological spaces indexed by $I$ is a collection of maps $\phi_r\colon X_r\rightarrow Y_r$ for each $r\in I$ such that the $\phi_r$ commute with the inclusions, i.e. $\phi_s \circ X_{r,s} = Y_{r,s}\circ\phi_r $ for all $r\leq s$.
\end{defn}

We will furthermore make use of the following notions for filtrations.
\begin{defn}
    Let $X$ be an $I$-filtered object (i.e.\ a vector space, simplicial complex, etc.).
    The \emph{birth} of an element $v\in X_r$ for some $r\in I$ is the infimum of all $s\in I$ such that $v\in X_s$.
    A point $p\in I$ is called a regular point for $X$, if there is a neighborhood $U$ of $p$ such that $X_{r,s}$ is an isomorphism for any $r,s\in U$. Otherwise, $p$ is called a critical point.
\end{defn}

The homology functor maps filtered topological spaces to persistence modules.

\begin{defn}
    A \emph{persistence module} indexed by $I\subset \R$ is a collection of vector spaces $ (V_r)_{r\in I} $ together with linear maps $V_{r,s}\colon V_r \rightarrow V_s$ for $r\leq s$ such that $V_{r,t} = V_{s,t}\circ V_{r,s}$ for $r\leq s\leq t$. The linear maps $V_{r,s}$ are called \emph{structure maps}.
\end{defn}

\begin{defn}
    Let $V$ and $W$ be persistence modules. A \emph{morphism of persistence modules} $\phi\colon V\rightarrow W$ is a collection of linear maps $\phi_r\colon V_r\rightarrow W_r$ that commute with the structure maps, i.e. 
    \[
    \phi_s \circ V_{r,s} = W_{r,s} \circ \phi_r,\quad \text{ for all } r,s\in I,\, r\leq s.
    \]
    It is an \emph{isomorphism of persistence modules} if all $\phi_r$ are isomorphisms.
    The \emph{image} $\im \phi$ of $\phi$ is the persistence module given by $(\im \phi_r)_{r\in I}$ with structure maps $(\im \phi)_{r,s} = W_{r,s}|_{\im {\phi_r}} $ being the appropriate restrictions of the $W_{r,s}$.
\end{defn}

We are now ready to give a more conceptual description of the cycling signature $\Cyc(\gamma,Y) = (\im H_1(i_{\gamma,r}))_{r\in I}$. By collecting the maps $ i_{\gamma,r} $ for $r\in I = [0,r_0(Y))$, we obtain a morphism 
\[
i_{\gamma}\colon O(\rho(\gamma)) \to \Delta Y
\]
of $I$-filtered topological spaces, where $\Delta Y=(\Delta Y_r)_{r\in I}, \Delta Y_r=Y$ for all $r\in I$, is the constant filtration of $Y$. Passing to homology induces a morphism 
\begin{equation}
    \label{eq:cycling-signature-as-image}
    H_1(i_{\gamma})\colon H_1(O(\rho(\gamma))) \to H_1(\Delta Y)
\end{equation}
of persistence modules. The cycling signature is the image of this morphism,
i.e., the morphism on the top row of the diagram

\begin{equation}
    \begin{tikzcd}
        H_1(O(\rho(\gamma))) \arrow[rr, "H_1(i_{\gamma})"]                    &  & H_1(\Delta Y)                             \\
                                                                              &  &                                           \\
        H_1(\Cech(\rho(\gamma))) \arrow[uu, "\simeq"] \arrow[rr, "H_1(\phi)"] &  & H_1(\Delta Y_\infty) \arrow[uu, "\simeq"]
    \end{tikzcd}
\label{eq:discretization-diagram}
\end{equation}

The bottom row of this diagram is a discretization of the top row and will be introduced in the following. We use \v{C}ech complexes to discretize thickenings of segments (\cref{subsec:internal-external-thickenings}), cubical complexes $Y_\infty$ to construct comparison spaces (\cref{subsec:cubical-comparison-space}) and acyclic carriers to obtain the chain map on the bottom of this diagram (\cref{subsec:cycling-map-discretization}).
The problem of computing a basis for a cycling signature then reduces to computing the image of a map of filtered vector spaces which is treated in \cref{subsec:maps-of-filtered-vector-spaces}.

\subsection{Background: \v{C}ech filtration and a nerve theorem}

We start by recalling some notions from algebraic topology, for a general introduction see \cite{Munkres84}.
Given a simplicial complex $K$, we denote its geometric realization by $|K|$. Note that `geometric realization' is functorial, i.e.\ it also applies to simplicial maps. We will furthermore make use of the barycentric subdivision $\sd K$ of a simplicial complex $K$. 

A standard method to construct a simplicial complex from a point cloud is the \v{C}ech complex.
\begin{defn}
    The \emph{\v{C}ech complex} $\Cech_r(P)$ of a finite set of points $P\subset M$ at radius $r$ in a metric space $(M,d)$ is the abstract simplicial complex
    \[
    \Cech_r(P) = \left\{ Q\subset P\mid \bigcap_{q\in Q} \overline{B(q,r)} \neq \emptyset \right\},
    \]
    where $B(q,r)$ is the ball with center $q$ and radius $r$ in that space. 
    We sometimes write $\Cech_r(P, M, d)$ to highlight the metric space. 
    The \emph{\v{C}ech filtration} of $P$ is the filtered simplicial complex $\Cech(P) = (\Cech_r(P))_{r\in [0,\infty)}$.
\end{defn}

Given a finite set of points $P\subset\R^d$, there are now two filtrations associated to it: the filtration of metric thickenings $O(P)$ and the \v{C}ech filtration $\Cech(P)$. In the following we state a nerve theorem from \cite{roll2023-unified} that guarantees topological equivalence of these two filtrations. We will furthermore use the explicit description of the homotopy equivalence for computations. 

\begin{figure}[H]
    \centering
    \includegraphics[width=0.65\linewidth]{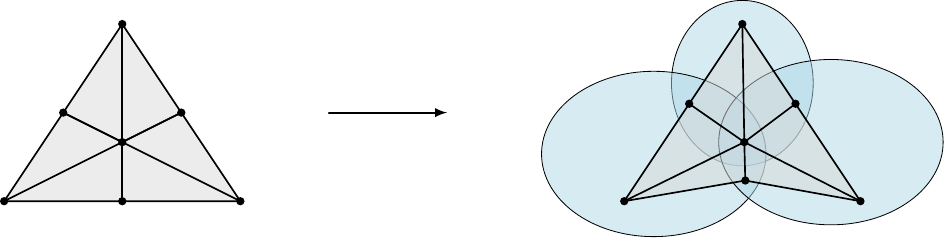}
    \caption{Illustration of the map $f_r$ (cf.\ \cite{roll2023-unified}).}
    \label{fig:nerve-map-illustration}
\end{figure}

In the following, we denote the filtration of the barycentric subdivisions of the \v{C}ech filtration by $K = (K_r)_{r\in[0,\infty)} = (\sd \Cech_r(P))_{r\in[0,\infty)}$. We now define a map $f\colon|K| \rightarrow O(P)$ from this filtration to the filtered topological space $O(P)$.
We start by defining $f$ on the 0-skeleton of $K$, i.e., the filtration of 0-skeleta. 
For this, fix any vertex $v$ in $\bigcup_r K_r$. 
Since $v$ is in the barycentric subdivision of $\Cech_r(P)$, it corresponds to a simplex $\sigma = \{ x_0,\dots, x_n\}$ spanned by vertices $ x_0,\dots, x_n$ in $P$. 
In particular, the birth radius of $v$ in $K$ is equal to the birth radius of $\sigma$ in $\Cech(P)$ which, by definition of the \v{C}ech complex, 
is the smallest radius $r_b>0$ such that 
\[
    \bigcap_{i=0}^n \overline{B(x_i,r_b)} \neq \emptyset.
\]
Fix any point $x$ in this intersection. We set $f_r(v) = x$ for all $r\geq r_b$, see \cref{fig:nerve-map-illustration}.
After repeating this for the entire 0-skeleton, we extend $f_r$ via linear interpolation to all of $|K_r|$, i.e., for any simplex $\sigma = \{x_0,\dots,x_n\}$ in $K_r$, we set, in barycentric coordinates,
\begin{equation}
    f_r(\lambda_0 x_0 + \dots + \lambda_n x_n) = \lambda_0 f(x_0) + \dots + \lambda_n f(x_n).
    \label{eq:nerve-map-defn}
\end{equation}

\begin{thm}
    For any point cloud $P\subset \R^d$, the map $f \colon |K|\rightarrow O(P)$ constructed above is a morphism of filtered topological spaces such that each $f_r$ is a homotopy equivalence.
    \label{thm:cech-nerve}
\end{thm}
\begin{proof}
    For any fixed $r$, the map $f_r$ is precisely the map $\Gamma$ from \cite{roll2023-unified} Theorem 3.1, therefore it is a homotopy equivalence.
    Furthermore, by Theorem 3.11 in the same paper, this homotopy equivalence is natural with respect to morphisms of pointed covered spaces. By Example 1.6 therein, these include the setting of barycentric subdivisions of \v{C}ech complexes.
\end{proof}

\begin{cor}
    For any finite point cloud $P\subset \R^d$, there is an isomorphism $ H_1(\Cech(P))\cong H_1(O(P)) $ of persistence modules.
\end{cor}

\subsection{Internal and External Thickenings}
\label{subsec:internal-external-thickenings}

We now treat the left vertical arrow of the diagram \eqref{eq:discretization-diagram}.
Thus far, the thickenings of segments $\gamma$ were taken internally in the unit tangent bundle, i.e.\ we considered
\[
    O_r^{I,C}(\rho(\gamma)) := O_r(\rho(\gamma)) = \{ p\in \UTB(X)\mid d_C(p,\rho(\gamma))\leq r  \},
\]
where $I$ indicates that the thickening is internal and $C$ is the metric parameter.
To apply \cref{thm:cech-nerve}, we want to work with external thickenings, where the balls are taken in $\T(\R^d) = \R^{2d} $, i.e.
\[
    O_r^{E,D}(\rho(\gamma)) = \{ p\in \T(X)\mid d_D(p,\rho(\gamma))\leq r  \}.
\]
As it turns out, with a suitable rescaling of the metric, these two thickenings are equivalent (see \cref{sec:appendix-internal-external-utb} for details).

\begin{prop}
    Let $\gamma\subset \R^d$ be finite. With $C \geq 0$, $r\in [0,C)$ and $D_r = r/(C\theta(r/C))$, the map
    \[
        \alpha_r \colon O_r^{E,D_r}(\rho(\gamma))\rightarrow O_r^{I,C}(\rho(\gamma)), \quad (p,v) \mapsto (p,v/\norm{v}_2)
    \]
    is a homotopy equivalence. The collection of maps $\alpha$ is furthermore a map of filtered topological spaces.
    \label{prop:main-external-internal-equiv}
\end{prop}
\begin{proof}
    See \cref{prop:appendix-external-internal-equiv}.
\end{proof}

\begin{cor}
    Let $\gamma\subset \R^d$ finite, $C \geq 0$ and $D_r = r/(C\theta(r/C))$. 
    There is a commutative diagram of $[0,C)$-indexed filtered topological spaces
    \[
        \begin{tikzcd}
        {O^{I,C}(\rho(\gamma))}                                    &  & {O^{E,D_\bullet}(\rho(\gamma))} \arrow[ll, "\alpha"'] \\
                                                                   &  &                                                       \\
        {|\Cech(\gamma,\UTB,d_C)|} \arrow[rr, "h"] \arrow[uu, "g"] &  & {|\Cech(\gamma,\T,d_{D_\bullet})|} \arrow[uu, "f"']  
        \end{tikzcd}
    \]
    where for fixed filtration value all maps are homotopy equivalences. We write $D_\bullet$ to emphasize that the parameter depends on $r$.
    \label{cor:cycling-i-e-equivalence}
\end{cor}
\begin{proof}
    Note that $\Cech_r(\gamma,\UTB,d_C) = \Cech_r(\gamma,\T,d_{D_r})$ for $r\in [0,C)$, since intersections of internal balls with radius $r$ and metric $d_C$ are nonempty if and only if external balls with radius $r$ and metric $d_{D(r)}$ are nonempty (see \cref{cor:i-e-intersection-equivalent}). 
    Thus, we can choose the same geometric realization for both of them.

    Since $O^{E,D_\bullet}(\rho(\gamma))$ is a family of subsets of $\R^{2d}$, by \cref{thm:cech-nerve}, the right vertical map can be chosen to be a map $f = (f_r)_r$ of filtered topological spaces where each component is a homotopy equivalence. 
    For the map $\alpha$, this follows from \cref{prop:main-external-internal-equiv}.

    We can thus construct the left vertical map as $g_r = \alpha_r\circ f_r\circ \operatorname{id}_r$.
\end{proof}

Applying homology to the map $g$ from the corollary, we get an isomorpism for the left vertical arrow in diagram \eqref{eq:discretization-diagram}.

\subsection{Cubical Comparison Spaces}
\label{subsec:cubical-comparison-space}

We now treat the right vertical arrow in diagram \eqref{eq:discretization-diagram} and explain how we construct the comparison space $Y$.
For this, we introduce the unit tangent bundle with respect to the $l_\infty$-norm
\[
    \UTB_\infty(\R^d) = \{ (p,v)\in \T(\R^d) : \norm{v}_\infty = 1 \}
\]
and note that there is a homeomorphism $\eta\colon \UTB_2(\R^d) \rightarrow \UTB_\infty(\R^d), (p,v)\mapsto (p,v/\norm{v}_\infty)$. We can thus obtain a comparison space $Y$ by finding a neighborhood $U$ of $\eta(\rho(\Gamma))$ and then taking $Y = \eta^{-1}(U)$. Computationally, this is convenient, since spheres $S_\infty^d$ are much easier to cover with boxes than spheres $S_2^d$.

Given a point $p\in \R^d$ and a radius $r\in \R_{>0}$, we write $Q_r(p) = \prod_{i=1}^d \left[p_i-\tfrac{r}{2},p_i+\tfrac{r}{2}\right]$ for the box with side lengths $r$ centered at $p$. 
To construct the comparison space, we use boxes of the form $Q_{r,k}(p,q) = Q_r(p)\times Q_{1/k}(q) \subset T_\infty(\R^d)$ where $ r,k\in \R_{>0} $, $p\in r\Z^d$  and $ q\in 1/k\Z^d $. 

Note that the collection of these boxes covers $ \R^{2d} $. 
To cover the unit tangent bundle we fix two parameters: the space box size $r \in \R_{>0}$ and the sphere box size $1/k$, where $k\in \N_{>0}$, and consider the following cover of $\UTB_\infty(\R^d)$
\[
\mathcal{Q}_{r,k}^d = \{ Q_{r,k}(p,q)\mid p\in r\Z^d,\ q\in 1/k\Z^d,\ \norm{q}_\infty = 1 \}.
\]

We can obtain comparison spaces by covering a neighborhood of $\eta(\rho(\Gamma))$, e.g. 
\[
    Y_\infty=  \{ Q\in \mathcal{Q}_{r,k}^d\mid Q\cap \eta(O_\epsilon(\rho(\Gamma)))\neq \emptyset \},\quad \text{for some } \epsilon \geq 0,
\]
and then taking $ \eta^{-1}( p_\infty(|Y_\infty|) )$, where $p_\infty\colon \T(\R^d) \setminus (\R^d\times \{ 0 \}) \rightarrow \UTB_\infty(\R^d), (p,v)\mapsto (p,v/\norm{v}_\infty)$ and we write $| Y_\infty |= \bigcup Y_\infty$.

\begin{defn}
    A \emph{cubical comparison space} is a set $Y_\infty \subset \mathcal{Q}_{r,k}^d $ for some $r>0$ and $k\in \N$ such that $Y = \eta^{-1}(p_\infty(|Y_\infty|))$ is a comparison space for $\Gamma\subset \R^d$.
    \label{defn:cubical-comparison-space}
\end{defn}

\begin{prop}
    For any $K \subset \mathcal{Q}_{r,k}^d$, where $r>0$ and $k\in \N_{>0}$ we have a homotopy equivalence $ p_\infty(|K|)\simeq |K| $.
    \label{prop:p-infty-hty-equivalence}
\end{prop}
\begin{proof}
 See \cref{prop:appendix-p-infty-hty-equivalence}
\end{proof}

Any $Q\in \mathcal{Q}_{r,k}^d$ is a polytope and thus has faces of all dimensions $\leq d$, we denote the collection of all faces by $\operatorname{\downarrow}Q$.
Given $K \subset \mathcal{Q}_{r,k}^d$, we write $\operatorname{\downarrow} K = \bigcup_{Q\in K} \operatorname{\downarrow}Q$. 

\subsection{Discretization of the Cycling Signature}
\label{subsec:cycling-map-discretization}
In this section, we introduce the map on the bottom row of the diagram \eqref{eq:discretization-diagram}.

Throughout this section, let $\Gamma\subset \R^d$ be a time series, $\gamma $ a segment of $\Gamma$ and $Y_\infty $ a cubical comparison space for $\Gamma$ with respect to $d_C$ and admissible radii $I\subset [0,C)$. We furthermore denote by $R\in I$ the largest critical point of $(\Cech_r(\rho(\gamma)))_{r\in I}$ in $I$.

\begin{prop}\label{prop:subdivision-to-carrier}
    Consider the map 
    \[
        \Theta \colon |\sd \Cech(\rho(\gamma))| \rightarrow \Delta Y_\infty  , \quad \Theta = \eta \circ \alpha\circ f,
    \]
    where $\alpha\circ f$ is as in \cref{cor:cycling-i-e-equivalence} and $\eta$ as in \cref{defn:cubical-comparison-space}.
    There is a subdivision $K_R$ of $\sd \Cech_R(\rho(\gamma))$ such that 
    \[
        F(\sigma) = \operatorname{\downarrow} \{ Q\in Y_\infty \mid Q \cap \Theta_{R}(\sigma) \neq \emptyset \}
    \]
    is an acyclic carrier (see \cite{Munkres84}, Theorem 13.3) from $K_R$ to $Y_\infty$.
\end{prop}
\begin{proof}
    \emph{Step 1:} We claim that there is a subdivision $K_R$ of $\sd \Cech_R(\rho(\gamma))$ with the property that for any face $H$ of the $l_\infty$ unit ball  $B_\infty(0,1) \subset \R^d $, the set $ \Theta_R^{-1}( \R^d\times H )$ is a subcomplex of $K_R$.

    The set $\gamma_{E}^{R} := O_R^{E,D_R}(\rho(\gamma)) $ is compact and does not intersect the zero section $\R^d\times \{ 0 \} \subset \T(X)$ since, with respect to $d_C$, the distance between $\R^d\times \{ 0 \}$ and $\UTB(\R^d)$ is exactly $C$ while $R < C$ by assumption. Thus there are constants $r_1,r_2,r_3$ such that 
    \[
        \gamma_{E}^{R} \subset \overline{B_\infty(0,r_1)}\times (\overline{B_\infty(0,r_2)}\setminus B_\infty(0,r_3)) =: A\subset \T(X).
    \]
    We claim that $A$ is the geometric realization of a simplicial complex $L$. In fact, the set $A$ can be decomposed into sets $A_H = A\cap (\R^d\times \{ rH\mid r>0 \} )$ where $H $ is a face of $ B_\infty(0,1) \subset \R^d $. 
    Each of the $A_H$ is a convex polyhedron since, pictorially, it is the product of the cube $B_\infty(0,r_1)$ with the pyramid $\{ rH\mid r>0 \}$ that has the corner at the origin chopped off. Furthermore, the intersection of two such sets is either empty, or a face of either of the polyhedra. Therefore, $A$ has a natural cell structure that can be subdivided into a simplicial complex $L$ such that $A = |L|$, see e.g. Prop. 2.9 in \cite{Rourke1972}.
    
    By \eqref{eq:nerve-map-defn}, $f_{R} \colon |\sd \Cech_{R}(\rho(\gamma))| \rightarrow \gamma_{E}^{R}\subset |L| $ is a linear map on each simplex. 
    Thus, there are simplicial subdivisions $K_{R}$ of $\sd \Cech_R(\rho(\gamma))$  and $L'$ of $L$ such that $f_{R} $ is a simplicial map with respect to these simplicial complexes (Lemma 2.13 in \cite{Rourke1972}).
    We turn $K_R$ into a filtered simplicial complex in the obvious way, i.e. by  $ K_r = \{ \sigma \in K_{R}\mid |\sigma| \subset |\Cech_r(\rho(\gamma)) | \}$.
    
    \emph{Step 2:} We need to show that $F(\sigma)$ is nonempty and acyclic, and that $F(\hat \sigma)\subset F(\sigma)$ if $\hat{\sigma}$ is a face of $\sigma$. Since the latter property is clear by definition, we now focus on the former.

    Let $\sigma \in K_{R}$ and set $\tau = \Theta_{R}(\sigma)$.  
    By above subdivision, we have $\tau\subset \R^d \times H $ for a facet $ H$ of $ B_\infty(0,1)$. 
    Therefore, restricted to $\sigma $, the map $\Theta_{R}$ is given by the composition of the linear map $f_{R}$ with $\eta\circ \alpha_{R}$, which in the tangent component is the radial projection onto $H$. 
    Both of these maps preserve compactness and convexity, thus $\tau$ is compact and convex. 
    It thus follows from \cref{cor:cover-by-utb-boxes-is-acyclic} that $F(\sigma)$ is acyclic.    
\end{proof}

\begin{cor}
    There is a commutative diagram
    \[
\begin{tikzcd}
H_1(|\sd \Cech(\rho(\gamma))|) \arrow[rr, "H_1(\Theta)"] &  & H_1(\Delta |Y_\infty|) \arrow[dd, "H_1(\psi_2)"'] \\
                                                         &  &                                                   \\
H_1(K) \arrow[uu, "H_1(\psi_1)"] \arrow[rr, "H_1(\nu)"]  &  & H_1(\Delta Y_\infty)                             
\end{tikzcd}
    \]
    where $K$ is a subdivision of $\sd \Cech(\rho(\gamma))$ and $\psi_1$ and $\psi_2$ are chain maps that induce the isomorphism between singular and cellular homology.
\end{cor}
\begin{proof}
    Let $K$ be the subdivision and $F$ be the acyclic carrier from Proposition \ref{prop:subdivision-to-carrier}. We get an algebraic acyclic carrier by setting $G(\sigma) = C(F(\sigma))$ to be the cellular chain complex of $F(\sigma)$.
    By the acyclic carrier theorem, there is a chain map $\phi_R \colon C(K_R) \rightarrow C(Y_\infty)$ carried by $G$. We turn this into a filtered map by the appropriate restrictions.

    Let $C(\Theta_R)$ denote the chain map induced by $\Theta_R$. We now claim that the chain map $\psi_2\circ C(\Theta_R) \circ\psi_1$ is also carried by $G$. To see this, recall that the isomorphism between singular and cellular homology has the property that for any subcomplex $ L_R $ of $K_R$, and any subcomplex $Z$ of $Y_\infty$ we have
    $
        \psi_1(L_R) \subset \psi_1(|L_R|),\, \psi_2(|Z|) \subset \psi_2(Z).
    $
    By construction, we have $\Theta_R(|\sigma|)\subset |F(\sigma)|$ for any simplex $\sigma \in K_R$ and thus
    \[
    \psi_2\circ C(\Theta_R)\circ\psi_1(C(\sigma)) \subset\psi_2\circ C(\Theta_R)(C(|\sigma|))\subset \psi_2(\Theta_R(|\sigma|))\subset \psi_2(|F(\sigma)|)\subset C(F(\sigma)).
    \]
    Thus, $\psi_2\circ C(\Theta_R) \circ\psi_1$ is chain homotopic to $\phi_R$. We can turn $\psi_2\circ C(\Theta_R) \circ\psi_1$ into a filtered chain map using the appropriate restrictions to obtain the top map of the diagram. Furthermore, we note that the chain homotopy between $\psi_2\circ C(\Theta_R) \circ\psi_1$ and $\phi_R$ respects the acyclic carrier, which implies that appropriate restrictions of those give a chain homotopy between $\psi_2\circ C(\Theta_r) \circ\psi_1$ and $\phi_r$ for every $r$. Therefore, the above diagram is a commutative diagram of persistence modules.
\end{proof}

\begin{cor}
    Let $Y = \eta^{-1}(p_\infty(|Y_\infty|))$. There is a commutative diagram
    \[
\begin{tikzcd}
H_1(O(\rho(\gamma))) \arrow[rr, "H_1(i_{\gamma})"]                    &  & H_1(\Delta Y) \arrow[dd, "\simeq"'] \\
                                                                      &  &                                     \\
H_1(\Cech(\rho(\gamma))) \arrow[uu, "\simeq"] \arrow[rr, "H_1(\phi)"] &  & H_1(\Delta Y_\infty)               
\end{tikzcd}
            \]
\label{cor:cycling-diagram}
\end{cor}
\begin{proof}
    Each space in this diagram is isomorphic to the corresponding space in the preceding corollary. 
    For the top left and top right spaces, this follows from \cref{thm:cech-nerve} and \cref{prop:p-infty-hty-equivalence}.
    For the bottom row, this follows from topological invariance of barycentric subdivision.
\end{proof}

\subsection{Maps of Filtered Vector Spaces}
\label{subsec:maps-of-filtered-vector-spaces}

We now explain how to compute a basis for the image of a map from a persistence module into a filtered vector space.
Throughout this section, all persistence modules are pointwise finite dimensional (i.e.\ all $V_i$ are finite dimensional) and indexed by $T = \{ 0,1,\dots, q\}$.

\begin{defn}[c.f.\ \cite{SchmahlLifespan2023}]
	Let $ V $ be a persistence module. Its \emph{immortal part} is the persistence module $ V^\infty =(V_i^\infty)_{i\in T}$ where 
	\[
	V_i^\infty = \{ v \in V_i \mid v = 0 \text{ or } V_{i,q}(v) \neq 0 \}, 
	\]
	its \emph{mortal part} is the persistence module $ V^\dagger=(V_i^\dagger)_{i\in T} $ where
	\[
	V_i^\dagger = \{ v \in V_i \mid V_{i,q}(v) = 0 \}, 
	\]
	and the structure maps are the appropriate restrictions. 
\end{defn}

A filtered vector space is equal to its immortal part and has trivial mortal part. More generally, the immortal part of any persistence module $V$ is isomorphic to a filtered vector space.

The cycling signature is expressed as the image of a morphism form a persistence module into a filtered vector space. In the following, we show that, computationally, it is enough to consider the case of a map between filtered vector spaces. 

\begin{prop}
	Let $ \phi \colon V \rightarrow W $ be a morphism of persistence modules. If $ W $ is a filtered vector space, then 
	\begin{enumerate}
		\item $ \im \phi $ is a filtered vector space, and
		\item $ \im \phi = \im \phi|_{V^\infty}$.
	\end{enumerate}
	\label{prop:image-morphism}
\end{prop}
\begin{proof}
	For the first claim, we need to show $ \im \phi_i \subset \im \phi_j $ for $ i \leq j $. Note that post-composition with an inclusion $W_{i,j}$ does not change the image, thus $\im \phi_i = \im W_{i,j} \circ\phi_i $. The claim now follows from
	\[
	\im \phi_i = \im W_{i,j} \circ\phi_i = \im \phi_j \circ V_{i,j} \subset \im \phi_j.
	\]
	
	For the second claim, note that a persistence module is a direct sum of its mortal and immortal part. It is therefore enough to show that $ \im \phi_{|V^\dagger} = 0 $. For this, fix $i \in T$ and note that for an arbitrary $ v \in V_i^\dagger $
	\[
	\phi_i(v) = W_{i,q}\circ\phi_i(v) = \phi_q \circ V_{i,q}(v) = 0
	\]
	where we again use that $ W_{i,q} $ is an inclusion map.
\end{proof}

As a consequence of the above proposition, the image of a morphism $\phi\colon V\rightarrow W$ of a persistence module into a filtered vector space can be computed by considering the restriction $\phi|_{V^\infty}$ which is a map of filtered vector spaces.
In the setting of vector spaces, the image of a linear map can be computed using matrix reduction as described in the following. Recall that $\F$ is some fixed field.

\begin{defn}
	The pivot $\operatorname{pivot}(v)\in\N$ of a vector $ v\in \F^m $ is the largest $ i\in 1:m$ such that $v_i \neq 0$ or zero if no such $i$ exists. A matrix is called \emph{reduced} if all non-zero pivots of its columns are distinct.
\end{defn}

Any matrix $M\in \F^{m\times n}$ can be reduced, i.e.\ there is an invertible and upper triangular matrix $ U\in \operatorname{Gl}_n(\F) $ such that $ R = M U \in \F^{m\times n}$ is reduced (see e.g.\ \cite{edelsbrunner2008}). 
Suppose we have a decomposition $R = M U$. Since $R$ is reduced, a basis for $\im M$ is given by $\{ R_i\mid \operatorname{pivot} R_i \neq 0 \}$. A different basis can be obtained from the columns of $M$. 

\begin{lem}
    \label{lem:matrix-reduction-image-basis}
    Suppose $R = M U$ with $M\in \F^{m\times n}$, $R\in \F^{m\times n}$ reduced and $U\in \operatorname{Gl}_n(\F)$ upper triangular. Then a basis of $\im M$ is given by $\{ M_i \mid \operatorname{pivot}(R_i) \neq 0 \}.$
\end{lem}
\begin{proof}
	Let $R_{1:i}$ and $ M_{1:i} $ be the matrices consisting of the first $i$ columns of $R$ and $M$, respectively. Since $U$ is invertible and upper triangular, $\im R_{1:i} = \im M_{1:i}$ for all~$i$. In particular, a column $R_{i+1}$ is linearly independent from $\im R_{1:i}$ if and only if a column $M_{i+1}$ is linearly independent from $\im M_{1:i}$. The claim now follows, since $\{ R_i\mid \operatorname{pivot}(R_i) \neq 0 \}$ is a basis for $\im M$.
\end{proof}

Given a linear map $\phi\colon V \rightarrow W$ of (unfiltered) vector spaces, a basis $v_1,\dots, v_n$ of $V$ and a basis $\alpha_1,\dots, \alpha_m$ of the dual space $W^*$ of $W$, a basis of $\im \phi$ can be computed by reducing the matrix $ M = [\alpha_i(\phi(v_j))]_{i,j} $. If $R = MU$ as above, a basis for $\im \phi$ is given by
\[
    \{ \phi(v_j) \mid \operatorname{pivot}(R_j) \neq 0 \},
\]
since the map $w\mapsto (\alpha_1(w),\dots,\alpha_m(w))$ induces an isomorphism $W\cong \F^m$.

This algorithm readily generalizes to the setting of a map $\phi:V\to W$ of filtered vector spaces by reducing the matrix of $\phi_q\colon V_q \rightarrow W_q$, i.e.\ at the final index. The only added requirement is that the basis elements need to be ordered by birth.

\begin{prop}
\label{prop:basis}
    Let $\phi\colon V\rightarrow W$ be a map of filtered vector spaces. Let furthermore $B^V$ and $B^{W^*}$ be bases for $V$ and the dual space $W^*$ of $W$, respectively, and suppose
    \begin{enumerate}
        \item $B_q^V = \{  v_1,\dots, v_n \}$, with $\birth(v_i)\leq \birth(v_j)$ for $i\leq j$,
        \item $B_q^{W^*} = \{ \alpha_1,\dots, \alpha_m \}$, with $\birth(\alpha_i)\leq \birth(\alpha_j)$ for $i\leq j$.
    \end{enumerate}
    Let $M = [\alpha_i(\phi_q(v_j))]_{i,j}$ and suppose $R = MU$ with $R$ reduced and $U$ invertible and upper triangular. Then, a basis $B$ for $\im \phi$ is given by
    \[
    B_i = \{ \phi_i(v_j) \mid j\leq i \text{ and } R_j \neq 0 \}.
	\]
\end{prop}
\begin{proof}
    Clearly, $B_i\subset B_j$ for $i\leq j$. Now fix $i\in T$ and let $R_{1:i}$ and $ M_{1:i} $ be the matrices consisting of the first $i$ columns of $R$ and $M$, respectively.   
    Since the map $w\mapsto (\alpha_1(w),\dots,\alpha_m(w))$ induces an isomorphism $W\cong \F^m$, it is enough to show that 
    $\tilde{B}_i = \{ M_j \mid j\leq i \text{ and } R_j \neq 0 \}$ is a basis of $\im M_{1:i}$. This follows from \cref{lem:matrix-reduction-image-basis} since we have $R_{1:i} = M_{1:i} U_{1:i,1:i}$.
\end{proof}

We can now use this result in order to compute the cycling signature in the following way: Let $\Gamma$ be a time series, $Y$ a comparison space for $\Gamma$ and $\gamma$ a segment of $\Gamma$.
Consider the \v{C}ech filtration $\Cech(\gamma)$. Since $ \gamma$ is a finite point cloud, there is a finite number of critical points $ 0 = r_0\leq r_1\leq r_2\leq \dots $ where the filtration changes. Let $q$ be the largest integer with $r_q < r_0(Y)$ and define a filtered vector space via
\[
    V_k = \left(H_1(\Cech(\gamma))^\infty\right)_{r_k},  
\]
i.e. by evaluating the filtered vector space $H_1(\Cech(\gamma))^\infty$ at the critical point $r_k$.
We then get a map $\phi \colon V \rightarrow \Delta H_1(Y)$ of filtered vector spaces indexed over $T = \{0,1,\dots,q\}$ to which we can apply \cref{prop:basis} to get a basis $\tilde{B}$ of $\im \phi=\Cyc(\gamma,Y)$. We then set $B_r = \tilde{B}_i$, where $i$ is the largest integer such that $r_i\leq r$. This yields a basis $B$ of $\Cyc_r(\gamma, Y)$.

\section{Algorithm and Implementation}

To compute a cycling signature $\Cyc(\gamma, Y_\infty) $ using \cref{prop:basis}, we need to compute and reduce the matrix $M_{i,j} = [ \alpha_i(\phi(v_j)) ]_{i,j}$ where
\begin{itemize}
    \item the $\alpha_i$ are a basis of $H^1(Y_\infty)$ for some cubical comparison space $Y_\infty$ of $\Gamma$, and
    \item the $v_j$ are a basis for $H_1(\Cech(\rho(\gamma))) $.
\end{itemize}
In the following, we explain how we these steps are implemented in \cite{Hien2024}.
The computation furthermore requires a choice of metric $d_C$ and knowledge of the maximal filtration radius $r_0$ which determines the interval of admissible radii $I=[0,r_0)$. 
In \cref{subsec:heuristics-C-r0} we provide heuristics for choosing these two values.

\subsection{Computation of Cubical Comparison Spaces and their Cohomology}
In practice, we construct the cubical comparison space $Y_\infty$ by fixing parameters $r,k\in\N$ and collecting all boxes in $\mathcal{Q}_{r,k}$ that contain a point from $\Gamma$. 
While this does not guarantee that the map $\Theta_r $ from \cref{prop:subdivision-to-carrier}, and thus also the chain map $\phi $ from \cref{cor:cycling-diagram},
is well-defined for $r>0$, we have implemented computational heuristics that enable us to construct the necessary map in homology (see \cref{sec:map-heuristics}).
We then compute a basis $\alpha_1,\ldots,\alpha_m$ of  $H^1(Y_\infty)$ using a version of the coreduction algorithm \cite{Mrozek2008}.
Since the computation of the basis only depends on $\Gamma$, but not on any specific segment, we can do this once and reuse it for all segments in $\Gamma$.

\subsection{Computing the Persistent Homology of Segments}
Instead of the \v{C}ech filtration, it is common practice in computational topology, to compute using the Vietoris-Rips filtration
\[
    \VR_r(P) = \{ \sigma \subset P \mid \operatorname{diam} \sigma \leq r \}.
\]
We adopt this heuristic in the following way.
First, a standard persistence code (in our case Ripserer \cite{Ripserer2020}, a Julia implementation of Ripser \cite{Bauer2021}) is used to compute a barcode $(b_i,d_i)_i$ and generators $G = (g_i)_i $ for $H_1(\VR(\rho(\gamma)))$.
Second, we make use of the fact that the 1-skeleta of \v{C}ech and Vietoris-Rips complexes are identical. 
Therefore, any chain in $G$
is also a chain in $\C_1(\Cech(\rho(\gamma))$.
We can thus regard the set $G$ as an approximate  generating set for $H_1(\Cech(\rho(\gamma))$.
We note that this is a heuristic, since in general the map $\C_1(\VR(\rho(\gamma))) \rightarrow \C_1(\Cech(\rho(\gamma) )$ does not induce a well-defined map in homology.

We then let $G_\infty$ to be all generators in $G$ that correspond to bars $(b_i,d_i]$ that contain the maximal filtration radius $r_0$, i.e.\ with $b_i < r_0 < d_i$, to obtain an 'approximate basis´ for the immortal part of  $H_1(\Cech(\rho(\gamma))$.

\subsection{Computation of the Map in Homology}
We now explain how we compute $\phi(c)$ for some chain $c\in C_1(\Cech(\rho(\gamma))$  where $\phi$ is the map from \cref{cor:cycling-diagram}.
By definition, $c = \sum_k\lambda_k e_k$, $\lambda_k\in\F$, for some edges $e_k$ in $\Cech(\rho(\gamma))$. 
By linearity, it is thus sufficient to compute the images $\phi(e_k)$ for all $k$.

So, fix some edge $e \in \Cech(\rho(\gamma)) $.
Formally, the map $\phi $ is the composition of a subdivision and the map induced by the acyclic carrier from \cref{prop:subdivision-to-carrier}.
Computationally, we subdivide $e$ such that two consecutive vertices of the subdivision lie in adjacent boxes and each segment lies on the same face of $\partial B_\infty(0,1)$ in the unit tangent component.
We note that this takes care of the subdivision that is necessary because of the nerve theorem \ref{thm:cech-nerve}.
We then construct a 1-chain as follows: each vertex in the subdivision is mapped to a distinguished point of a box in $Y_\infty$ that contains it.
Each edge $e = [v,w]$ in the subdivision is then mapped to a chain that is contained in the union of the boxes and has $\phi(w) - \phi(v)$ as its boundary (this is how a chain map is induced by an acyclic carrier, see \cite{Munkres84}).
The sum of all these chains is then the image of the edge in $\C_1(Y_\infty)$.

\subsection{Computation of the Image Barcode}
Both, the previously computed basis $\alpha_1,\dots, \alpha_m$ of $H^1(Y_\infty)$ and the chains $\phi(c_1),\dots, \phi(c_n)$, where $G_\infty = \{  c_1,\dots, c_n\}$, are represented as vectors in $\F^k$ where $k$ is the number of edges in $Y_\infty$. The entries of $M_{i,j} = [ \alpha_i(\phi(v_j)) ]_{i,j}$ are then computed simply by taking the inner product of the previously computed vectors corresponding to $\alpha_i$ and $\phi(v_j)$.
In our examples, the matrix dimensions of $M$ are less than $5$, thus we reduce $M$ using a basic implementation of the standard reduction algorithm and then read off a basis of the cycling signature as explained in \cref{prop:basis} and the paragraph following it.

\subsection{Heuristics for the metric \texorpdfstring{$d_C$}{C} and the maximal filtration radius \texorpdfstring{$r_0$}{r0}}
\label{subsec:heuristics-C-r0}
To compute a cycling signature, we need to specify values for the constant $C$ of $d_C$ and the maximal filtration radius $r_0(Y)$. 
Given a cubical comparison space $Y_\infty \subset \mathcal{Q}_{r,k}$, we suggest $C = rk$ and $r_0 = r$ as default values. We now give a heuristic explanation for these choices.
Our choices are motivated by \cref{prop:thickening-of-cubical-sets} which states that a cubical complex consisting of radius $R$ cubes can be thickened by $ r < R$ without changing its homotopy type. 

Suppose $Y_\infty$ is constructed as above, i.e.\ consists of all boxes in $\mathcal{Q}^d_{r,k}$ that intersect with $\rho(\Gamma)$. 
Similar to \cref{prop:thickening-of-cubical-sets}, we can consider
\[
    Y_\infty^{s,C} = \{ O_s(Q,d_C) \mid Q \in Y_\infty \}
\]
where all the individual cubes are thickened with respect to $d_C$.
Recall, that the cubes in $\mathcal{Q}^d_{r,k}$ are of the form $ Q_r(p)\times Q_{1/k}(q) $. 
Thus, to have a thickening that preserves the nerve, the first component should be thickened less than $r$, the second component should be thickened by less than $1/k$. 
For the two components of the metric $d_C$ from \cref{eq:UT-metric}, this yields the equations $ s < r $ and $ Cs < 1/k $. 
Thus, with $C = rk$, we have that any thickening with $ s < r $ preserves the nerve of $Y_\infty^{s,C}$. 
Furthermore, any other thickening that has this guarantee is contained in $ |Y_\infty^{s,C}| $, in this sense our choices are therefore optimal.

\subsection{Heuristics for the Map in Homology}
\label{sec:map-heuristics}
As described in \cref{sec:computation}, the construction of the chain selector does not take into account the thickened boxes but only works with the original box cover.
It thus may happen in the computation of the map in homology that a segment $e$ of an edge in $\Cech(\rho(\gamma))$ is contained in a box that is not contained in $Y_\infty$. 
By out choice of parameters, this segment is contained in the union of the thickened boxes. 
As a heuristic, we then use the collection of these boxes in the chain map construction.

\section{Discussion}

Cyclic signatures provide a novel framework in which to identify coarse clusters of cycling motions and transitions between these coarse clusters from time series data generated by an ODE.
Furthermore, we claim that there is a natural association with our construction and fundamental structure of dynamics of ODEs.
Fixed points are trivial examples of recurrent dynamics.
The flow box theorem states that a point  in phase space $x\in X$ is either a fixed point or that in a neighborhood $N$ of $x$ the dynamics is diffeomorphic to parallel flow, i.e., there is no local recurrence.
Our insistence that the comparison space be a subset of a unit tangent bundle naturally enforces the condition that cycling signatures do not identify trivial recurrence.
Again, citing the flow box theorem, non-trivial recurrent dynamics can only arise if some point in the neighborhood $N$ eventually carried back by the flow to the neighborhood $N$.
By definition a nontrivial cycling signature captures this phenomenon.
Of course, this does not imply that all recurrence will be captured by cycling signatures; the Hopf fibration provides an example for which a sufficiently long noisy trajectory will result in only trivial cycling signatures.
Though we do not have a precise statement about the extent of applicability of cycling signatures, the examples of this paper show that it can provide useful information and a heuristic as to its applicability is as follows.

The cycling signature is an invariant that can be assigned to any time series segment. 
Statistics of these signatures provide a qualitative and quantitative description of cycling motions in the data.
The resulting analysis is global in nature; clusters of cycling motions are detectable and transitions between cycling motions can be analyzed. 

Our approach reduces the identification of cycling behavior to linear algebra; cycling signatures are linear subspaces of $H_1(Y)$.  While a canonical basis would facilitate a more direct correspondence between cycling signatures and cycling motions, this does not seem to be possible.  Furthermore, the linear dependence among cycling motions (cf.~\cref{fig:s21-inclusion}) suggests that this might not be desirable.
For example, the dimension of $H_1(Y)$ in the Dadras system is 5, while we observe 6 different cycling classes which are all contained in a 4-dimensional subspace of $H_1(Y)$.
With this in mind, our pipeline was developed to analyze cycling motions without the need for a canonical basis.

In this work, we define the cycling signature as the image of a map in degree 1 homology and make use of 1d and 2d cycling spaces to analyze cycling dynamics in the example systems.
We do believe, however, that higher rank signatures can contain useful information. Investigating this is topic of further work.
We furthermore did not investigate the natural generalization of cycling to higher homology groups in this work even though this may be necessary to analyze high-dimensional time series data.

Cycling signatures depend on several inputs: the time series $\Gamma$, the choice of comparison space $Y$ and the metric $d_C$ on the unit tangent bundle.
A rigorous understanding of how cycling signatures change as function of these inputs exceeds the scope of this paper.
From an empirical perspective, cycling signatures appear to be robust (see \cref{sec:lorenz-stability-experiments}).

A key ingredient in computing $H_1(Y)$ is the construction of the comparison space~$Y$. Here, we constructed $Y$ from an extrinsic cubical covering of a thickening of the given time series.  The advantage of this is that typically the number of cubes is much smaller than the number of points in the given time series so that the computation of $H_1(Y)$ becomes feasible. 
Furthermore, this construction is applicable to general time series data and provides a certain amount of flexibility through the box size parameters.
Of course, other constructions are possible. 
Potentially, a different construction of $Y$ could lead to a simplified computational pipeline. 
Furthermore, for some applications of our technique it may be necessary to design other types of comparison spaces that are adapted to the system and data of interest.

Our new tool is particularly useful for higher dimensional systems for which visualizations are of limited use. 
At first glance it appears challenging to compute cycling signatures in a higher dimensional system (the complexity of homology computations depends on the number of elementary cells whose number in general scales exponentially in the dimension of state space).
However, cycling signatures only require knowledge of homology in degree 1.
Therefore,  we only need to consider 0-, 1- and 2-dimensional cells. 
What is more, our approach is flexible with respect to the type of cell complex we use.
We thus believe that, by a suitable improvement of our computational pipeline, the treatment of systems of considerably higher dimension should be possible.

\section*{Acknowledgments}

This research has been supported by the German Research Foundation (DFG – Deutsche Forschungsgemeinschaft TRR109 Project ID 195170736). K.M. was partially supported by the National Science Foundation under awards DMS-1839294 and HDR TRIPODS award CCF-1934924, DARPA contract HR0011-16-2-0033,  National Institutes of Health award R01 GM126555, and Air Force Office of Scientific Research under award number FA9550-23-1-0011 and FA9550-23-1-0400. K.M. was also supported by a grant from the Simons Foundation.

\bibliographystyle{abbrv} 
\bibliography{cycling-paper}

\clearpage

\appendix

\section{Generation of the time series data}
\label{sec:data-generation}
The computational results in \cref{sec:application} were obtained using our Julia implementation   \texttt{CyclingSignatures.jl} \cite{Hien2024}, differential equations were integrated using the Julia package  \texttt{DifferentialEquations.jl} \cite{rackauckas2017differentialequations}. 

\subsection{Lorenz System}
\label{sec:Lorenz}

The classical Lorenz system \cite{Lo-63} is a  system of differential equations which exhibits a characteristic chaotic attractor. The equations are
\begin{align*}
	\dot{x} &= \sigma(y-x),\\
	\dot{y} &= x(\rho-z)-y,\\
	\dot{z} &= xy-\beta z.
\end{align*}
and we consider the system at the classical values $ \sigma=10 $, $\rho=28$ and $ \beta = \frac{8}{3} $.
The time series was generated by numerically integrating the system using the Tsit5 scheme \cite{TSITOURAS2011770} with initial condition $(0,10,0)$ and step size control such that two subsequent time series points have distance less than $1$.

\subsection{Perturbed Double Well}

We consider a stochastic differential equation
\begin{equation}
	dx = f(x) dt + \sigma dW,
	\label{eq:SDE}
\end{equation}
$x=(q,p)\in\R^2$, with $W$ denoting the Wiener process.  The vector field $f:\R^2\to\R^2$ is an asymmetric double well Hamiltonian system  with two homoclinic orbits, given by
\begin{align*}
	f(x) &= \begin{pmatrix}H_p(x)\\ -H_q(x)\end{pmatrix} - a h(H(x)) \begin{pmatrix}H_q(x)\\ H_p(x)\end{pmatrix}
\end{align*}
where $h(z) = (z^3-z)/2$ and 
\[
H(q,p) = p^2/2 + q^4/8 - q^2/2 - q^3/15 - q/10.
\]
For the experiments in \cref{sec:application}, we used $a=0.02$, for the time series in \cref{fig:dw-cycling-comparison}, we set $a=0$.
The time series were generated by integrating \eqref{eq:SDE} using the SRIW1 scheme \cite{Roessler2010} with the initial value $x=(1,0.75)$, step size 0.01 and noise level $\sigma=0.015$. 

\subsection{Dadras System}

In \cite{Dadras.2012}, the authors consider a 4d system of differential equations with a single equilibrium at the origin and a 4 wing hyperchaotic attractor. The equations are
\begin{align*}
	\dot{x} &= ax-yz+w,\\
	\dot{y} &= xz-by,\\
	\dot{z} &= xy-cz+xw\\
	\dot{w} &= -y,
\end{align*}
where we set $ (a,b,c)=(8,40,14.9)  $ as suggested in the original article. In order to confine the system to a neighborhood of the origin, we use the nonlinear rescaling $x \mapsto x/\sqrt{\norm{x}_2}$ for all time series points. 
The time series was generated by numerically integrating the equation using the Tsit5 scheme \cite{TSITOURAS2011770} with initial condition $(10,1,10,1)$ and step size control such that two subsequent time series points have distance less than $ 0.8 $.

\section{Lorenz Stability Experiments}
\label{sec:lorenz-stability-experiments}

The computation of cycling signatures depends on the choice of comparison space and unit tangent bundle metric.  
From a practical point of view it is important that the computation of cycling signatures is stable with respect to these choices.
Lacking a thorough theoretical understanding of the stability of cycling signatures, we experimentally investigate the effect of small parameter changes on the presented results. 

\begin{figure}[h]
	\centering
	\includegraphics[width=\linewidth]{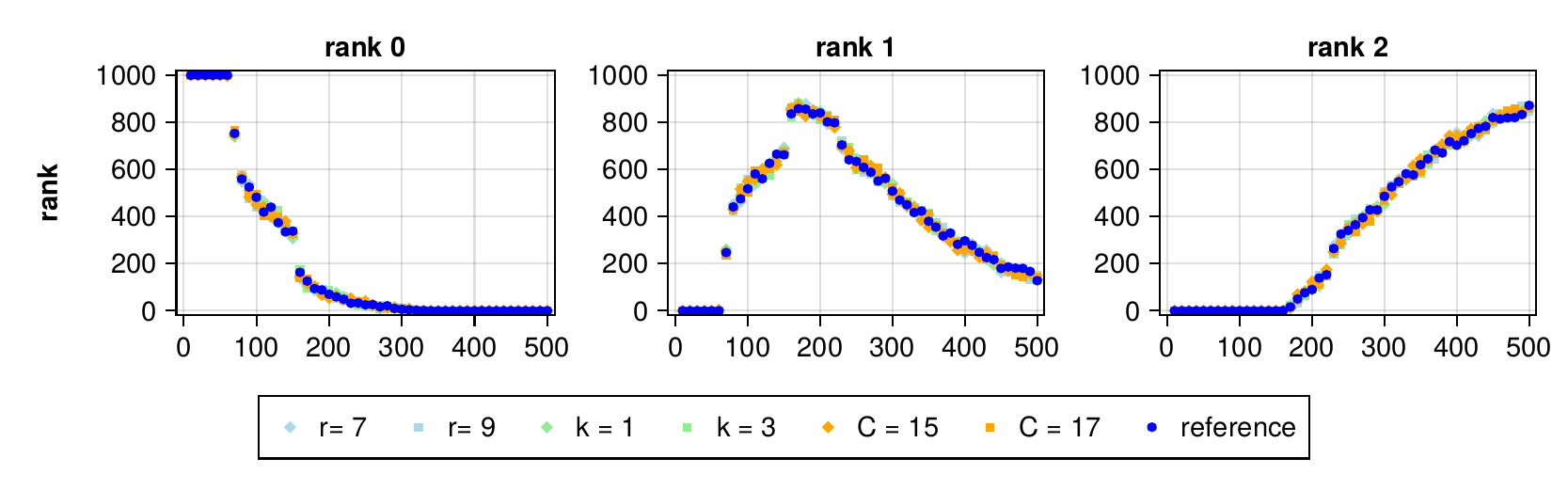}
	\caption{Rank distribution versus segment time span for a variety of parameter choices. All cycling signatures are evaluated at $r=4$.}
	\label{fig:stability-rank}
\end{figure}

\begin{figure}[h]
	\centering
	\includegraphics[width=\linewidth]{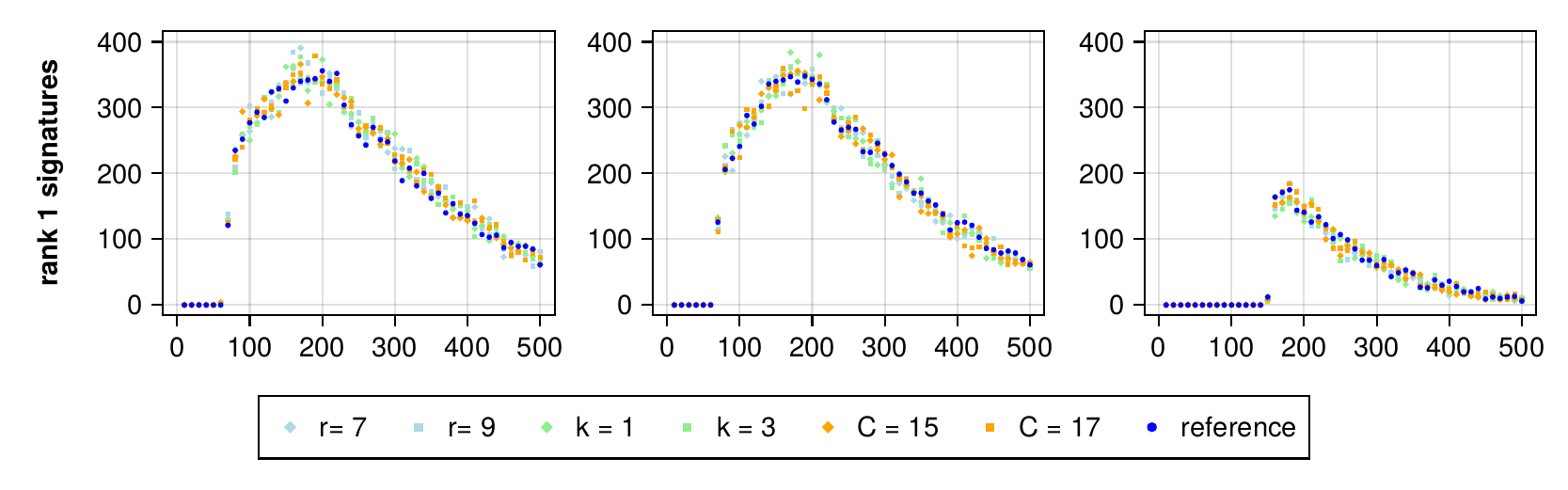}
	\caption{Rank 1 signature distribution versus segment time span for a variety of parameter choices. All cycling signatures are evaluated at $r=4$.}
	\label{fig:stability-sig1}
\end{figure}

We consider the Lorenz time series that we already investigated in \cref{sec:application}. 
There, we chose the cover parameter $(r,k)$ to be $(8,2)$ and the metric constant $C$ to be $16$.
To demonstrate stability of our results, we repeat the experiment with different parameters.
More precisely, we vary the parameters $r,k,C$ as follows:
\begin{itemize}
    \item $r=7,\,8,\,9$ with $k = 1$ and $C = 16$,
    \item $k=1,\,2,\,3$ with $r = 8 $ and $C = 16$,
	\item $C = 16,\,17,\,18 $ with $r=8$ and $k=2$.
\end{itemize}
For each parameter combination, we sample $1000$ segments for segment time spans $\{10, 20, \dots, 500\}$.
Thus, changes in the cycling signatures arise not only from changes in parameters but also due to sampling.
As we can see in \cref{fig:stability-rank} and \cref{fig:stability-sig1}, changes in the statistics of rank distribution and are small in the tested range of parameters. 

\section{Covers of Convex Compact Sets}
\label{sec:appendix-covers}

We consider the set
\[
\mathcal{Q}^k(\R^d) = \left\{ \prod_{i=1}^d [v_i,v_i+ w_i] \mid v\in \Z^d,\, w\in \{0,1\}^d,\text{ such that } \sum_i w_i = k  \right\} 
\]
of all $k$-dimensional unit cubes in $\R^d$ that have vertices in the integer lattice.
We furthermore set $\mathcal{Q}(\R^d) = \bigcup_{k=0}^d \mathcal{Q}^k(\R^d)$.
Given $P,Q \in \mathcal{Q}(\R^d)$, we call $P$ a face of $Q$ if $P\subset Q$ and write $\operatorname{\downarrow} Q $ for the set of all faces of $Q$.
For $K\subset \mathcal{Q}(\R^d)$ we set $\operatorname{\downarrow} K = \bigcup_{Q\in K}\operatorname{\downarrow} Q$ and write $|K| = \bigcup_{Q\in K} Q$. 
A set $K\subset \mathcal{Q}(\R^d)$ is called a \emph{cubical (grid) complex} if $K = \operatorname{\downarrow} K$.
The outer cubical cover of a set $X\subset \R^d$ is the cubical complex
\[
     K(X) = \operatorname{\downarrow}\left\{ Q \in \mathcal{Q}^d(\R^d) \mid Q\cap X \neq \emptyset \right\}.
\]

\begin{prop}
    The outer cubical cover $K(\sigma)$ of a convex compact set $\sigma$ is acyclic.
    \label{prop:cover-of-simplex-is-acyclic}
\end{prop}

Before the proof, we state two corollaries. 
For this, let $A\in \operatorname{GL}(d)$ and set $\mathcal{Q}^k_A(\R^d) = \{ A Q \mid Q\in \mathcal{Q}^k \}$ and $ \mathcal{Q}_A(\R^d) = \bigcup_{k=0}^d \mathcal{Q}^k_A(\R^d) $ where $A Q = \{ Aq\mid q\in Q\}$.

\begin{defn}
    The \emph{outer cubical cover of $X\subset \R^d$ with respect to $\mathcal{Q}_A(\R^d)$} is the cubical complex
    \[
         K_A(X) = \operatorname{\downarrow}\left\{ Q\in \mathcal{Q}_A^d(\R^d) \mid Q \cap X \neq \emptyset \right\}.
    \]
    \label{def:outer-cubical-cover}
\end{defn}

\begin{cor}
    Let $A \in \operatorname{GL}(d)$ and suppose $\sigma$ is convex and compact. 
    Then $K_A(\sigma)$ is acyclic.
    \label{cor:cover-of-simplex-acyclic}
\end{cor}
\begin{proof}
    Linear automorphisms preserve compactness, convexity and homology.
    By \cref{prop:cover-of-simplex-is-acyclic}, the outer cubical cover $K(A^{-1} \sigma )$ is acyclic. The claim now follows from $|K_A(\sigma)| = A |K(A^{-1}\sigma)| $. 
\end{proof}

The following is the corollary that is used in the construction of the acyclic carrier in \cref{prop:subdivision-to-carrier}. 
Recall the notation $\mathcal{Q}_{r,k}^d$ from \cref{subsec:cubical-comparison-space}. 

\begin{cor}
    Let $\sigma\subset \UTB_\infty(\R^d)$ be compact and convex. For any $r\in \R,k\in \Z$, the set
    \[
        L = \operatorname{\downarrow} \{ Q\in \mathcal{Q}_{r,k}^d \mid Q\cap \sigma \neq \emptyset \}
    \]
    is acyclic.
    \label{cor:cover-by-utb-boxes-is-acyclic}
\end{cor}
\begin{proof}
    Let $v\in \R^{2d}$ be the vector with first $d$ entries $r$ followed by $d$ entries $\tfrac{1}{k}$. 
    Let furthermore $A\in \R^{2d\times 2d}$ be the diagonal matrix with $v$ on the diagonal.
    Define $\tau = \{ w - v/2 \mid w\in \sigma \}$.
    By \cref{cor:cover-of-simplex-acyclic}, the outer cubical cover $ K_A(\tau) $ is acyclic. 
    
    We now claim that $ |L| $ and $ |K_A(\tau)| $ differ only by a translation by $v/2$.
    For this, fix $Q\in K_A(\tau) $ and set  $\tilde{Q} = \{ w+v/2\mid w\in Q \}$.
    Clearly, $\tilde{Q}\cap \sigma \neq \emptyset$.
    We thus only need to show that $\tilde{Q}$ is contained in $\mathcal{Q}_{r,k}^d$.
    It is easy to see that $ \tilde{Q} = Q_r(p)\times Q_{1/k}(q)$ for some $p\in r\Z^d$ and $q\in 1/k\Z^d$, thus it is only left to show that $\norm{q}_\infty = 1$.
    To see this, it is enough to show that whenever $\norm{q}_\infty \neq 1$, we have $Q_{1/k}(q) \cap S_\infty^1 = \emptyset$.
    Pictorially, the box $Q_{1/k}(q)$ is a box of side length $1/k$ centered at $q\in 1/k\Z^d$, thus we have $ q - 1/(2k)\leq \norm{x}_\infty \leq q + 1/(2k)$ for all $x\in Q_{1/k}(q)$. In particular, $\norm{x}_\infty\neq 1$ if $q\in 1/k\Z^d$ with $\norm{q}_\infty\neq 1$.

    Since translation is a homeomorphism, it now follows that $ |L| $ is acyclic.
\end{proof}

The proof of \cref{prop:cover-of-simplex-is-acyclic} will proceed by induction on the number of maximal cubes in $K(\sigma)$. In the induction step, a layer of cubes will be cut off via a Mayer-Vietoris argument. 

We need several lemmas in preparation for proving the proposition.

\begin{defn}
	A convex, compact set $\sigma\subset\R^d$ is in \emph{general position} with respect to $\mathcal{Q}^k(\R^d)$ if $\sigma \cap \Int Q \neq \emptyset $ for all $ Q\in K_d(\sigma)$.
\end{defn}

As an example, $\sigma = \{0\}\subset \R$ is not in general position, since it does not intersect the interior of the interval $[0,1]\in K(\sigma)$. As it turns out, any compact convex set can be enlarged slightly to another compact convex set that is in general position and has the same outer cover. We prove this in the next two lemmas.

\begin{lem}
	Let $\sigma \subset \R^d$ be convex and compact. Then there is $ \bar\epsilon > 0$ such that for all $\epsilon \in [0,\bar{\epsilon})$,
	\[
		K(\sigma) = K(\sigma^\epsilon),
	\]
\end{lem}
\begin{proof}
    Let $L = \operatorname{\downarrow}(K(\R^d)\setminus K(\sigma))$ be the cubical complex of all boxes that do not intersect $\sigma$.
	Any box $ Q \in L $ has positive distance from $\sigma$, i.e., it is $d_2(Q,\sigma) > 0$. Furthermore, given any $ C > 0 $, there are only finitely many boxes with $ d(Q,\sigma) < C $. In particular, there is a minimal distance from $\sigma$ to $L$, i.e. the minimum
	\[
		\bar\epsilon = \min \{\operatorname{dist}(\sigma,Q)\mid Q\in L \}
	\]
	exists. Let $ \epsilon\in [0,\bar{\epsilon}) $. By the triangle inequality, any box in $ L $ has positive distance to $ \sigma^\epsilon $, thus $ K(\sigma) = K(\sigma^\epsilon) $.
\end{proof}

\begin{lem}
	Let $\sigma \subset \R^d$ be convex and compact. Then there is an $\epsilon > 0 $ such that
	\begin{enumerate}
		\item $K(\sigma) = K(\sigma^\epsilon)$, and
		\item $\sigma^\epsilon$ is in general position.
	\end{enumerate}
    \label{lem:general-position}
\end{lem}
\begin{proof}
	The previous lemma implies the existence of a $\bar{\epsilon}$ such that that the first condition is satisfied for all $ \epsilon \in [0,\bar{\epsilon})$. We now fix any such $\epsilon$.

	For the second condition, let $ Q \in K(\sigma) $ be arbitrary. We pick any point $x \in \sigma \cap Q \neq \emptyset $ and note that $ \emptyset \neq B(x,\epsilon) \cap \Int Q \subset \sigma^\epsilon \cap \Int Q $. Thus, $\sigma^\epsilon$ is in general position.
\end{proof}

\begin{proof}[Proof of \cref{prop:cover-of-simplex-is-acyclic}]
    By \cref{lem:general-position}, we can assume that $ \sigma$ is in general position.
    
    We prove the claim by induction on the number of top-dimensional cubes in $ K(\sigma) $. 
    The case $ n=1 $ is simple, since a single cube is contractible.
    Now suppose $K(\sigma)$ contains $ n+1$ top-dimensional cubes.
    Let $d$ denote the dimension of the top-dimensional cubes.
    Since $n>1$, there must be a $ k \in \N$ such that projection on the $k$-th coordinate satisfies $ p_k(K(\sigma)) = [m_1,m_2]$ with $m_2 - m_1 \geq 2$. Now let $ l = m_2-1 $. 
    
    We will now `cut' $K(\sigma)$ at the hyperplane $H = p_k^{-1}(\{l\})$. Setting $A = \operatorname{\downarrow}\{ Q\in K_d(\sigma)\mid p_k(Q)\subset (-\infty,l] \} $ and $B = \operatorname{\downarrow}\{ Q\in K_d(\sigma)\mid p_k(Q)\subset [l,\infty) \} $, we have the short exact sequence of chain complexes
	\[
		0 \rightarrow C(A \cap B) \rightarrow C(A)\oplus C(B) \rightarrow C(A \cup B) \rightarrow 0
	\]
	where the last term is equal to $ C(K(\sigma)) $. Its long exact sequence is the Mayer-Vietoris sequence in homology. It thus suffices to show that the left and middle term in this sequence are acyclic.

    We start by showing that $A$ and $ B$ are acyclic. In order to apply the induction hypothesis, we argue that $A$ and $B$ are outer covers of compact convex sets. 
    By general position, for any $ Q \in K_d(\sigma) $, we can fix a point $ x_Q \in \sigma \cap \Int Q $.
    Define 
	\[
		\delta_A = \frac{1}{2} \max \{ l - p_k(x_Q) \mid Q \in A \}, \text{ and } \delta_B = \frac{1}{2} \min \{ p_k(x_Q) - l \mid Q \in B \}.
	\]
    Since the $x_Q $ are contained in the interior of maximal cubes, we have $p_k(x_Q)\not\in \Z$ and thus $\delta_A,\delta_B> 0$.
    We set $\sigma_A = \sigma \cap p_k^{-1}((-\infty,l-\delta_A]) $ and $\sigma_B = \sigma \cap p_k^{-1}([l+\delta_B,\infty) $ and claim $A = K(\sigma_A)$ and $B = K(\sigma_B)$.

    To show $A\subset K(\sigma_A)$, fix a d-dimensional cube $Q\in A$. 
    We have previously fixed a point $x_Q\in \sigma \cap \operatorname{Int} Q$. 
    An elementary estimate shows $p_k(x_Q) < l-\delta_A$, thus we have $ x_Q \in \sigma_A $. 
    Since $Q$ is the unique cube with $x_Q$ in its interior, we must have $Q\in K(\sigma_A)$.
    Conversely, if $Q\in K(\sigma_A)$, we have $Q\in K(\sigma)$, by the definition of outer cubical covers.
    Since $\sigma_A \subset p_k^{-1}((-\infty,l-\delta_A])$, we must have $Q \subset p_k^{-1}((-\infty,l])$ and thus $Q\in A$.
    This shows that also $A\supset K(\sigma_A)$.
    Similar arguments work to show $B = K(\sigma_B)$.

	We now have that $ A $ and $B$ are outer covers of compact convex sets. As we show next, both sets have less than $n+1$ top-dimensional cubes. To see this, note that any $d$-dimensional cube $Q\in K(\sigma)$ is either contained in $A$ or in $B$ and that, by the definition of $l$, both sets contain at least one $d$-dimensional cube.
    By induction, $C(A)$ and $C(B)$ are therefore acyclic. 
 
    It is left to show that $ A \cap B $ is acyclic. 
    For this, note that $\hat{\sigma} = \sigma \cap H$ is compact and convex. Furthermore, $\hat{\sigma} \subset A\cap B$.
    Now let 
    \[
    \pi_k\colon \R^d\rightarrow\R^{d-1},\quad (x_1,x_2,\dots,x_d)\mapsto (x_1,\dots,\hat{x}_k,\dots, x_d)
    \]
    be the projection omitting the k-th component. Then $\pi(\hat{\sigma})$ is compact, convex, and covered by
    \[
        C = \{ \pi_k(Q)\mid Q\in A\cap B \}.
    \]
    We claim $C = K(\pi(\hat{\sigma}))$. This follows from convexity of $\sigma$ as follows: Fix $Q\in K(C)$. To it corresponds a unique $d-1$-dimensional elementary cube $P$ in $\R^d$ that is in $A\cap B$. 
    In particular there are $d$-dimensional elementary cubes $Q_A\in A$ and $Q_B\in B$ that have $P$ as their face. 
    This means that there are points $x_A\in Q_A\cap \sigma$ and $x_B\in Q_B\cap \sigma$. 
    By convexity of $\sigma$, the line segment between them lies in $\sigma$ and intersects $P$ at some point $x$. 
    This point satisfies $\pi_k(x)\in Q$.
    Thus the claim follows.
    
    Furthermore, since any $d-1$ cube of $K(C)$ is adjacent to a $d$-dimensional cube in each of $A$ and $B$, we have that $2\# K_{d-1}(C) \leq n $. Thus, by induction, $C$ is acyclic, and since $C$ and $A\cap B$ are homeomorphic, so is $A\cap B$. With the Mayer-Vietoris argument from above, this finishes the proof.
\end{proof}

\section{Cubical Comparison Space Induces Comparison Space}

Throughout this section, fix $d,k\in\N_{>0}$ and let $\pi_\infty\colon \R^d\setminus \{ 0\} \rightarrow S^1_\infty,\, x\mapsto x/\norm{x}_\infty$.
The goal of this appendix is to prove \cref{prop:p-infty-hty-equivalence}. 
The proof will make use of the following cover of the relevant cubes. 
For any cube $Q = Q_{1/k}(p)$, where $p\in \Z^d$ with $\norm{k}_\infty$, and any $i = 1,\dots,d$, we define
\[
    \operatorname{Dec}_i(Q) = \begin{cases}
            \{\{ x\in Q \mid |x_i|\geq |x_j| \text{ for } j\neq i \}\} & \text{if } |p_i| = k, \\
            \emptyset & \text{else}.
    \end{cases}
\]
and $\operatorname{Dec}(Q) = \bigcup_{i=1}^d \operatorname{Dec}_i(Q) $.

\begin{lem}
    For any cube $Q = Q_{1/k}(p)\subset \R^d$, where $p\in \Z^d$ with $\norm{p}_\infty=k$, we have $Q \subset |\operatorname{Dec}(Q)| $.
    \label{lem:appendix-dec-covers-cube}
\end{lem}
\begin{proof}
    Note that any $x\in Q$ can be written as $x = p/k + v$ with $\norm{v}_\infty = \frac{1}{2k}$. 
    Thus, if $|p_i|=k$, it is $|x_i|=|p_i/k+v_i|\geq 1-\frac{1}{2k}$. 
    Otherwise, $|p_i|\neq k$ and we have $|x_i|=|p_i/k+v_i|\leq (k-1)/k +\frac{1}{2k} = 1 - \frac{1}{2k}$.

    Thus, for any $x\in Q$, there is an $i$ such that $|x_i|\geq |x_j|$ for $j\neq i$ and $|p_i| = k$.
    For this index, we have $x \in | \operatorname{Dec}_i(Q) | $.
\end{proof}

\begin{defn}
    Let $H = \{ x\in \R^d \mid x_i = a \}$ with $a \neq 0$ be a hyperplane.
    The \emph{radial projection (from the origin)} to $ H $ is the map
    \[
        p\colon \R^d \setminus \{ x\in \R^d \mid x_i = 0 \} \rightarrow H, \quad x\mapsto a\frac{x}{x_i}.
    \]
    \label{lem:radial-projection-of-compact-set-is-compact}
\end{defn}

\begin{lem}
    Let $H = \{ x\in \R^d \mid x_i = a \}$ with $a \neq 0$ be a hyperplane and $p$ the radial projection from the origin onto $H$.
    Let $K$ be a convex set with $ K \cap H = \emptyset$.
    Then $p(K)$ is convex.
\end{lem}
\begin{proof}
    Fix $x,y\in \R^d$ and note that for $\lambda\in [0,1]$
    \[
        p(\lambda x +(1-\lambda)y) = \frac{\lambda x_i}{\lambda x_i +(1-\lambda)y_i} \frac{ax}{x_i} + \frac{(1-\lambda) y_i}{\lambda x_i +(1-\lambda)y_i} \frac{ay}{y_i} =: a(\lambda) p(x) + b(\lambda) p(y),
    \]
    where we define $a(\lambda)$ and $b(\lambda)$. Note that $a$ and $b$ are well-defined, since either $x_i,y_i>0$ or $x_i,y_i<0$. 
    Furthermore, $a(\lambda)+b(\lambda)=1$ and $a(0) = b(1) = 0$, $a(1) = b(0) = 1$.
    The image of $\lambda \mapsto p(\lambda x +(1-\lambda)y)$ is thus the straight line between $p(x)$ and $p(y)$. From this observation, we see that $p(K)$ is convex.
\end{proof}

\begin{lem}
    For any cube $Q = Q_{1/k}(p)\subset \R^d$, where $p\in \Z^d$ with $\norm{p}_\infty=k$, and any $i=1,\dots,d,$ such that $\operatorname{Dec}_i(Q)\neq\emptyset$, it is
    \begin{enumerate}
        \item $R_i \in \operatorname{Dec}_i(Q) $ is compact and convex, 
        \item $\pi_\infty(R_i) \subset \{ x\in \R^d \mid x_i = p_i/|p_i| \}$, and
        \item $\pi_\infty(R_i)$ is compact convex.
    \end{enumerate}
    \label{lem:appendix-pi-infty-convexity}
\end{lem}
\begin{proof}
    \emph{(1):} $R_i$ is compact since it is the intersection of a compact cube with a (closed) halfspace. 
    Convexity follows from the observation that $r_i$ is the intersection of the convex sets $Q$ and $\{ x\mid |x_i|\geq |x_j| \text{ for } j\neq i\}$.

    \emph{(2)}: Let $x = \pi(y)\in \pi_\infty(R_i)$. We need to show that $x_i = p_i/|p_i|$. 
    We have
    \[
        x = \pi_\infty(y) = \frac{y}{\norm{y}_\infty} = \frac{y}{|y_i|}.
    \]
    Note that $y = p/k + v$ for a $v \in [-1/2k,1/2k]^d $. Furthermore, since $\operatorname{Dec}_i(Q)\neq\emptyset$, we have $|p_i| = k$. 
    Therefore, $x_i = (p_i + v_i)/ |p_i+v_i| = p_i/|p_i|$.
    
    \emph{(3)}: By (2), the restriction of $p_\infty$ to $R_i$ is a radial projection onto a hyperplane. Thus, convexity is preserved and $R_i$ is compact as the image of a compact set.
\end{proof}

\begin{lem}
    Let $P \subset \{ p\in \Z^d \mid \norm{p}_\infty = k \}$. Then
    \begin{enumerate}
        \item $\bigcap_{p\in P} Q_{1/k}(p)\neq \emptyset\,$ if and only if $\, \norm{p-q}_\infty\leq 1$ for all $p,q\in P$,
        \item If $\bigcap_{p\in P} Q_{1/k}(p)\neq \emptyset\,$, then $\bigcap_{A\in\mathcal{A}} A \neq \emptyset$ for any $\mathcal{A} \subset \bigcup_{p\in P} \operatorname{Dec}(Q_{1/k}(p))$,
        \item For any $\mathcal{A} \subset \bigcup_{p\in P} \operatorname{Dec}(Q_{1/k}(p))$ it is
        \[
            \bigcap_{A\in\mathcal{A}} A = \emptyset \implies \bigcap_{A\in\mathcal{A}} p_\infty(A) = \emptyset.
        \]
    \end{enumerate}
    \label{lem:appendix-cube-projection-intersection-prpoerties}
\end{lem}
\begin{proof}
    \emph{(1):} This is a simple calculation. 
    
    \emph{(2):} Set $a_i = \min \{ p_i\mid p\in P\} $ and  $b_i = \max \{ p_i\mid p\in P\} $. By (1), either $a_i=b_i$ or $a_i = b_i-1$. 
    Define 
    \[
        x_i = \begin{cases}
            \frac{a_i+b_i}{2k} & \text{if } a_i \neq b_i \text{ and } |b_i|\neq k, \\
            \frac{2b_i - \operatorname{sign}(b_i)}{2k} & \text{else.}
        \end{cases}
    \]

    We claim that $x\in A$ for every $A\in \mathcal{A}$. To see this, note that $|x_i - p_i|\leq 1/(2k)$ for any $p\in P$. 
    Therefore $x\in Q_{1/k}(p)$.
    Furthermore, if $|p_i| = k$, it is $|x_i| = 1 - \frac{1}{2k} $ and if $|p_i| < k$ it is $|x_i| \leq 1 - \frac{1}{2k} $.
    Thus, $x\in R$ for any $R\in \operatorname{Dec}(Q_{1/k}(p))$. In particular, $x\in \bigcap_{A\in\mathcal{A}} A$.

    \emph{(3):} Suppose $\bigcap_{A\in\mathcal{A}} A = \emptyset$. 
    By (1) and (2), there are $A,B \in \bigcap_{A\in\mathcal{A}} A$ with $A\in \operatorname{Dec}(Q_{1/k}(p))$ and $A\in \operatorname{Dec}(Q_{1/k}(q))$
    such that $p,q\in P$ and $\norm{p-q}_\infty > 1$. 

    Fix $i$ such that $|p_i - q_i|\geq 2$. 
    We will show 
    \begin{equation}
        \max \{ a_i \mid a \in \pi_\infty(A)\} < \min \{ b_i \mid b\in \pi_\infty(B) \}
        \label{eq:complicated-inequality}
    \end{equation}
    which implies $\pi_\infty(A) \cap \pi_\infty(B) = \emptyset$.

    For $x\in Q_{1/k}(p)$ and $y = \pi_\infty(x) = x/\norm{x}_\infty$, elementary estimates yield 
    \[
        y_i \in \begin{cases}
            \left[\frac{2p_i - 1}{2k+1}, \frac{2p_i + 1}{2k-1} \right] & \text{if } p_i \geq 1,  \\
            \\[-10pt]
            \left[\frac{- 1}{2k-1}, \frac{1}{2k+1} \right] & \text{if } p_i  = 0,  \\
            \\[-10pt]
            \left[\frac{2p_i + 1}{2k-1}, \frac{2p_i - 1}{2k+1} \right] & \text{if } p_i \leq -1. \\
        \end{cases}
    \]
     
    Without loss of generality, we assume $p_i + 2 \leq q_i$. 
    After substituting the relevant interval bounds into \cref{eq:complicated-inequality}, we check that the inequality is satisfied for all relevant values of $p_i$ and $q_i$:
    See \href{https://www.wolframalpha.com/input?i=%282q-1%29%2F%282k%2B1%29+-+%282p-1%29%2F%282k-1%29+%3E+0%2C+with+k%3E%3D2%2C+k%3E%3Dq+%2C+p%3E%3D+1%2C+p%2B2%3C%3Dq}{here} if $p_i,q_i\geq 1$ (or, by symmetry, if $p_i,q_i\leq -1$), 
    \href{https://www.wolframalpha.com/input?i=solve+%282q-1%29%2F%282k%2B1%29+-1%2F%282k%2B1%29+%3E+0%2C+with+k%3E%3D1%2C+2%3C%3Dq%2C+q+%3C%3D+k+}{here} if $ p_i = 0$ (or, by symmetry, $q_i =0$).
    In case $p_i < 0$ and $0 < q_i$, the inequality can easily seen to be true by comparing signs.
\end{proof}

We are finally able to prove the relevant proposition.

\begin{prop}
    For any $K \subset \mathcal{Q}_{r,k}^d$, where $r>0$ and $k\in \N_{>0}$ we have a homotopy equivalence $ p_\infty(|K|)\simeq |K| $.
    \label{prop:appendix-p-infty-hty-equivalence}
\end{prop}
\begin{proof}
    Let $S = \{ (p,k) \mid Q_r(p)\times Q_{1/k}(q) \in K \}$.
    We let $\mathcal{C} = \{ Q_r(p)\times A \mid (p,q)\in S,\, A \in \operatorname{Dec}( Q_{1/k}(q) )  \}$.
    By \cref{lem:appendix-dec-covers-cube}, each cube in $K$ is covered by its corresponding elements in $\mathcal{C}$.
    Thus, the union of these cubes $|K|$ is covered by $\mathcal{C}$.
    This implies that $\mathcal{D} = \{ \pi_\infty(X)\mid X\in \mathcal{C} \}$ is a cover of $\pi_\infty(|K|)$.
    By \cref{lem:appendix-pi-infty-convexity}, the covers $\mathcal{C}$ and $\mathcal{D}$ consist of compact convex sets.

    We now claim that these two covers have the same Nerve. 
    For this, we need to show for any $X_i,\dots, X_n \in  \mathcal{C}$ 
    \[
        \bigcap_i X_i \neq \emptyset \iff \bigcap_i p_\infty(X_i) \neq \emptyset
    \]
    The implication `$\Rightarrow$' is obvious. The implication `$\Leftarrow$' follows from (3) in \cref{lem:appendix-cube-projection-intersection-prpoerties} together with the observation that the product of two cubes intersects if and only if the component cubes intersect.

    We want to apply the Nerve theorem from \cite{roll2023-unified}, Theorem 3.11 to above covers. For this, we turn the covers into pointed covers by choosing any point in $x_A$ for $A\in \mathcal{C}$ and the point $ p_\infty(x_A) $ for $p_\infty(A)$. 
    We then have that $p_\infty$ induces a morphism of pointed covers, each of which satisfies the assumptions of the Nerve theorem.
    This yields that $p_\infty$ is a homotopy equivalence between the covered spaces.
\end{proof}

\section{Thickenings of Cubical Complexes}
\begin{prop}
    Let $A$ be a positive multiple of the identity matrix, i.e. $A = rI_d$ for some $r>0$, and $K \subset \mathcal{Q}_{A}^d(\R^d)$. Set
    \[
        K^s = \{ O_s(Q) \mid Q\in K \}
    \]
    where $0 < s < r$ and the thickening is with respect to the Euclidean metric. Then $\bigcup_{Q\in K^s} Q \simeq | K |$.
    \label{prop:thickening-of-cubical-sets}
\end{prop}
\begin{proof}
    Note that two non-adjacent boxes in $\mathcal{Q}_{A}^d(\R^d)$ can only intersect if $s\geq r$. 
    Thus, for $s<r$, we have that $Q_0,\dots, Q_k \in K$ have nonempty intersection if and only if $ O_s(Q_0), \dots,  O_s(Q_k)\in K^s$ have nonempty intersection. Therefore, the nerves of the covers
    $K^s$ of $\bigcup_{Q\in K^s}Q $ and  $K$ of $ | K |$ are isomorphic. By the nerve theorem for closed convex sets, the spaces are therefore homotopy equivalent.
\end{proof}

\section{Internal and External Covers of Spheres}

Let $X\subset S_2^{d-1}$ be a point cloud on the $\ell^2$-unit sphere in $\R^d$. Any point $p\in S^{d-1}_2$ has an external (closed) ball
\[
    B(p,r) = \{ x\in \R^d\mid \norm{x-p}_2\leq r \}
\]
and an internal (closed) ball 
\[
    B^S(p,r) = \{ x\in S_2^{d-1}\mid \norm{x-p}_2\leq r\}.
\]

The goal of this section is to show that the external and internal thickenings of $X$ are equivalent, or more precisely, that we have 
\begin{equation}
    \bigcup_{p\in X} B^S(p,r) \simeq \bigcup_{p\in X} B(p,\theta(r)),\, \text{ where } \theta(r) = \sqrt{1-\left(1-\frac{r^2}{2}\right)^2 }.
    \label{eq:cover-hty-equiv-and-theta}
\end{equation}

\begin{prop}
    Let $0\leq r < 1$. The maps 
    \begin{equation}
        \pi\colon \bigcup_{p\in X} B(p,\theta(r)) \rightarrow \bigcup_{p\in X} B^S(p,r),\quad x\mapsto \frac{x}{\norm{x}_2}
        \label{eq:projection-ball-to-sphere}
    \end{equation}
    and
    \begin{equation}
        \psi\colon \bigcup_{p\in X} B^S(p,r) \rightarrow \bigcup_{p\in X} B(p,\theta(r)),\quad y\mapsto \frac{2-r^2}{2} y.
        \label{eq:map-sphere-to-ball}
    \end{equation}
    are homotopy equivalences and homotopy inverses to each other. 
    \label{prop:phi-psi-hty-equivalences}
\end{prop}

\begin{proof}
    We first consider the setting of a single ball, i.e., $X=\{p\}$ with $p\in S_2^{d-1}$.
    
    \textbf{Claim: $\pi$ is well-defined:} Fix $x\in B(p,\theta(r))$. Let $\lambda = \langle x,p\rangle/\norm{x}_2^2$ and note that $\langle x,p-\lambda x \rangle = 0$. 
    This implies
    \begin{align*}
        \norm{x -p}_2^2 &= \norm{(1-\lambda)x + \lambda x - p }_2^2 \\
        &= \norm{(1-\lambda)x}_2^2 - 2 \langle (1-\lambda)x,\lambda x - p\rangle+ \norm{\lambda x-p}_2^2\\
        &= \norm{(1-\lambda)x}_2^2 + \norm{\lambda x-p}_2^2 \geq\norm{\lambda x-p}_2^2.
    \end{align*}
    Using $\norm{p}_2=1$ and again $\langle x,p-\lambda x \rangle = 0$, we furthermore have
    \[
        \norm{\lambda x - p}_2^2 = \norm{\lambda x}_2^2 - 2\langle \lambda x, p\rangle + \norm{p}_2^2 =\norm{\lambda x}_2^2 -2\langle \lambda x, p-\lambda x + \lambda x \rangle +1 = 1 - \norm{\lambda x}_2^2.
    \]
    Combining the last two calculations, we get 
    \[
         \norm{\lambda x}_2^2 = 1 - \norm{\lambda x-p}_2^2 \geq 1 - \norm{x-p}_2^2\geq 1 - \theta(r)^2 = \left(1 - \frac{r}{2} \right)^2.
    \]
    We estimate the distance of $\pi(x)$ to $p$ using the preceding inequality
    \begin{align*}
        \norm{\pi(x) - p}_2^2 = \norm{\pi(x)}_2^2 - 2\langle \pi(x),p\rangle + \norm{p}_2^2  &= 1 - 2\frac{\langle x,p\rangle}{\norm{x}_2} + 1\\ 
        &= 2 - \lambda\norm{x}_2 
        \leq 2 - 2 \left(1 - \frac{r^2}{2} \right) = r^2
    \end{align*}
    which implies $\pi(x)\in B^S(p,r)$.

    \textbf{Claim: $\psi$ is well-defined.} Fix $y\in B^S(p,r)$. Note that
    \[
        \norm{y-p}_2^2 = \norm{y}_2^2-2\langle y,p \rangle+\norm{p}_2^2=2-2\langle y,p \rangle.
    \]    
    We calculate
    \begin{align*}
        \norm{\psi(y)-p}_2^2 &= \norm{\frac{2-r^2}{2}y-p}_2^2
        = \left(\frac{(2-r^2)}{2} \right)^2\norm{y}_2^2 - (2-r^2)\langle y,p\rangle + \langle p,p\rangle\\
        &= \frac{(2-r^2)^2}{4} - (2-r^2)\langle y,p\rangle + 1
        = \frac{(2-r^2)^2}{4} - (2-r^2)\frac{(2-\norm{y-p}_2^2)}{2} + 1\\
        &\leq  \frac{(2-r)^2}{4} - (2-r^2)\frac{2-r^2}{2} +1
        = 1 + \frac{(2-r^2)^2}{4} - \frac{(2-r^2)^2}{2} = 1 - \frac{(2-r^2)^2}{4}\\
        &= \theta(r)^2.\\
    \end{align*}

    \textbf{Claim: $\pi\circ \psi$ is the identity.} It is clear by definition that $\pi\circ \psi = \operatorname{id}$. 

    \textbf{Claim: $ \psi\circ\pi \cong \operatorname{id}$.}
    The straight-line homotopy $H\colon B(p,\theta(r))\times [0,1]\rightarrow B(p,\theta(r))$
    \[
        H(x,t) = (1-t)x + t\psi\circ\pi(x) =(1-t)x+t \frac{2-r^2}{2\norm{x}_2}x
    \]
    is well-defined since Euclidean balls are convex.

    \textbf{The general case: }
    All constructions so far extend to the setting of multiple balls. 
\end{proof}

\section{Internal and External Covers of the Unit Tangent Bundle}
\label{sec:appendix-internal-external-utb}
To construct the cycling signature, we consider metric thickenings of segments $\gamma$ in the unit tangent bundle $\UTB(\R^d)$. These thickenings are with respect to a metric
\[
    d_C((p,v),(q,w)) = \max\{\norm{p-q}_2, C\norm{v-w}_2\}.
\]

With an internal cover, we mean the sets $O_r(\gamma)^r$. This is comprised of \emph{internal balls} of the form
$
    \overline{B_C^I(x,r)} = \{ y\in \UTB(\R^d)\mid d_C(y,x)\leq r \}
$
and has the form
\[
    O_r(\gamma) = \bigcup_{x \in \gamma} \overline{B_C^I(x,r)}.
\]
The external cover consists of balls
$
    \overline{B_D^E(x,r)} = \{ y\in \T(\R^d)\mid d_D(y,x)\leq r \}
$
and has the form
\[
    \bigcup_{x \in \gamma} \overline{B_D^E(x,r)}.
\]

The goal of this section is to show that internal and external covers are equivalent. This will use the results from the previous appendix. Recall the definition of $\theta$ in \cref{eq:cover-hty-equiv-and-theta}.

\begin{prop}
    Let $C > 0$, $r\in [0,C)$ and $D = r/(C\theta(r/C))$. Then for $\gamma\in \UTB(\R^d)$ finite, the map
    \[
        \alpha\colon \bigcup_{x \in \gamma} \overline{B_D^E(x,r)} \rightarrow
        \bigcup_{x \in \gamma} \overline{B_C^I(x,r)},\quad (p,v) \mapsto (p,v/\norm{v}_2)
    \]
    is a homotopy equivalence.
    \label{prop:appendix-external-internal-equiv}
\end{prop}
\begin{proof}
    Recall the maps $\pi$ and $\psi$ from \cref{eq:projection-ball-to-sphere,eq:map-sphere-to-ball}. Using these, we have $\alpha((p,v)) = (p,\pi(v))$. Define furthermore
    \[
        \beta\colon\bigcup_{x \in \gamma} \overline{B_C^E(x,r)} \rightarrow
        \bigcup_{x \in \gamma} \overline{B_D^I(x,r)},\quad (p,v) \mapsto (p,\psi(v)).
    \]
    We now claim that $\alpha$ and $\beta$ are well-defined and homotopy inverse to each other.
    By definition of the $d_C$, we have
    \[
        \overline{B_D^E((p,v),r)} =\overline{B(p,r)}\times \overline{B(v,r/D)}= \overline{B(p,r)}\times \overline{B(v,\theta(r/C))}
    \]
    and similarly 
    \[
        \overline{B_C^I((p,v),r)} =\overline{B(p,r)}\times \overline{B^S(v,r/C)}.
    \]
    Thus, by \cref{prop:phi-psi-hty-equivalences}, both maps are well-defined. We furthermore have
    $\beta\circ\alpha(p,v) = (p,\pi\circ\psi(v)) = (p,v)$ and $ \psi\circ\pi \cong \operatorname{id} $ where a homotopy is given by the straight line homotopy in the second component, i.e.
    \[
        H((p,v)),t) = (p,(1-t)v + t\psi\circ\pi(v)).
    \]    
\end{proof}

\begin{cor}
    Let $C > 0$, $r\in [0,C)$ and $D = r/(C\theta(r/C))$. Then for $\gamma\in \UTB(\R^d)$ finite,
    \[
        \bigcap_{x \in \gamma} \overline{B_D^E(x,r)} \neq \emptyset \iff \bigcap_{x \in \gamma} \overline{B_C^I(x,r)} \neq \emptyset
    \]
    \label{cor:i-e-intersection-equivalent}
\end{cor}
\begin{proof}
    Fix $x\in \gamma$. It follows from the proposition and its proof that
    \[
        y\in \overline{B_D^E(x,r)} \implies \alpha(y)\in\overline{B_C^I(x,r)},\quad\text{ and }\quad y\in\overline{B_C^I(x,r)} \implies \beta(y)\in \overline{B_D^E(x,r)}.
    \]
    The equivalence now follows by mapping the intersection points using $\alpha$ and $\beta$.
\end{proof}

\end{document}